\documentclass[12pt]{article} 
\usepackage{amsmath}
\usepackage{amssymb} 
\usepackage{amscd}
\usepackage[pdftex]{graphicx}
\DeclareGraphicsExtensions{.pdf}
\usepackage[T1]{fontenc}
\usepackage{xspace}

\input psfig.sty

\newtheorem{theorem}{Theorem}

\newtheorem{proposition}[theorem]{Proposition}
\newtheorem{lemma}[theorem]{Lemma}
\newtheorem{corollary}[theorem]{Corollary}

\newtheorem{definition}[theorem]{Definition}

\def\id{{\mbox{Id}}}
\def\ch{{\mbox{ch}}}
\def\im{{\mbox{Im}}}

\def\Lie{{\mbox{Lie}}}

\def\cala{{\cal A}} 
\def\calb{{\cal B}} 
 
\def\call{{\cal L}}
\def\calm{{\cal M}} 
\def\caln{{\cal N}} 
\def\calo{{\cal O}} 

\def\cals{{\cal S}}
\def\calc{{\cal C}}

\def\frach{{\mathfrak h}}
\def\fracd{{\mathfrak d}} 
\def\fracS{{\mathfrak {S}}}
  \def\fracd{{\mathfrak d}} 
 \def\fracsu{{\mathfrak {su}}} 
 \def\fracsl{{\mathfrak {sl}}}
\def\fracso{{\mathfrak {so}}}
\def\calg{{C_{\mathrm{alg}}}}
\def\bbbone{\mbox{\rm 1\hspace {-.6em} l}}

\def\ug{{\mathbf u}}
\def\vg{{\mathbf v}}

\def\vg{{\mathbf v}}

\numberwithin{equation}{section}

\begin{document}
\enlargethispage{3cm}

 \thispagestyle{empty}
\begin{center}
{\large\bf MODULI SPACE AND STRUCTURE}
\end{center} 
 \begin{center}
{\large\bf OF NONCOMMUTATIVE 3-SPHERES}
\end{center} 
   
\vspace{0.3cm}

\begin{center} Alain CONNES \footnote{Coll\`ege de France, 3 rue
d'Ulm, 75  005 Paris, \\ I.H.E.S. and Vanderbilt University,
connes$@$ihes.fr}
and 
 Michel DUBOIS-VIOLETTE
\footnote{Laboratoire de Physique Th\'eorique, UMR 8627, Universit\'e Paris XI,
B\^atiment 210, F-91 405 Orsay Cedex, France, 
Michel.Dubois-Violette$@$th.u-psud.fr\\
}
\end{center} \vspace{0,5cm}

\begin{abstract}

 We analyse the moduli space and the structure of noncommutative 3-spheres.
We develop
 the notion of central quadratic form
for quadratic algebras, and 
 prove a general algebraic result 
which considerably refines the classical homomorphism 
from a quadratic algebra to a cross-product algebra
associated to the characteristic variety and lands
in a richer cross-product. 
It allows to control the $C^\ast$-norm
on involutive quadratic algebras 
and to construct the differential calculus in 
the desired generality.
The moduli space of noncommutative 3-spheres is identified with
equivalence classes of pairs of points in a symmetric space
of unitary unimodular symmetric  matrices. The  scaling foliation of the moduli space is identified to the gradient flow of
the character of a virtual representation of $SO(6)$. Its generic orbits are
connected components of real parts of elliptic curves 
which form a net of biquadratic curves with $8$ points in common. We show that generically these curves are
the same as the characteristic variety of the associated
quadratic algebra. We then apply the general theory 
of central quadratic forms to show that the noncommutative 3-spheres
 admit a natural ramified covering $\pi$ by a noncommutative
$3$-dimensional nilmanifold.
 This yields the differential calculus. We 
then compute the Jacobian of the ramified covering $\pi$
by pairing the direct image of the fundamental class of the
noncommutative
$3$--dimensional nilmanifold
with the Chern character of the defining unitary 
and obtain the answer as the product of a period (of
an  elliptic integral) by a rational function.
Finally we show that the hyperfinite 
factor of type II$_1$ appears as cross-product of the field $K_q$
of meromorphic functions on an elliptic curve by a subgroup of its Galois group
Aut$_{\mathbb C}(K_q)$.
 \end{abstract}

\vspace{0,5cm}

\noindent LPT-ORSAY 03-34 and IHES/M/03/56

\newpage
\baselineskip=0.7cm 
\tableofcontents

\section{Introduction}

\noindent Noncommutative Differential Geometry is a growing subject centering around
 the exploration of a new kind of geometric spaces which
do not belong to the classical geometric world. The  theory is resting 
on large classes of examples as well as on the elaboration of new general concepts.
The main sources of examples so far have been provided by:
\begin{itemize}

\item[1)] A general principle allowing to understand difficult quotients
of classical spaces, typically spaces of leaves of foliations, as 
noncommutative spaces. 

\item[2)] Deformation theory which provides rich sources of examples
in particular in the context of ``quantization" problems.

\item[3)] Spaces of direct relevance in physics such as the Brillouin zone
in the quantum Hall effect, or even space-time as in the context of 
 the standard model with its Higgs sector.
\end{itemize}

\noindent We recently came across a whole class of new noncommutative spaces
defined as  solutions of a basic equation of $K$-theoretic origin. 
This equation was at first expected to admit only  commutative 
(or nearly commutative) solutions. 

\noindent  Whereas  classical spheres provide simple solutions of arbitrary dimension $d$,
it turned out that
when the dimension $d$ is $\geq 3$ there are very interesting new,
and highly noncommutative, solutions. The first examples were given 
in \cite{ac-lan:2001}, and in \cite{ac-mdv:2002a} (hereafter refered to as Part I)
we began the classification of all solutions in the $3$-dimensional case, 
by giving an exhaustive list of noncommutative 3-spheres $S^3_\ug$, and analysing 
the ``critical" cases. 
We also explained in \cite{ac-mdv:2002a}
a basic relation, for generic values of our ``modulus" $\ug$, between the 
algebra of coordinates on the noncommutative $4$-space 
of which $S^3_\ug$ is the unit sphere, and the Sklyanin algebras
which were introduced  in the context of totally integrable systems.\\

\noindent In this work we analyse the structure of noncommutative $3$-spheres
$S^3_\ug$ and of their moduli space. 
We started from the above relation with the  Sklyanin algebras
and first computed basic cyclic cohomology invariants using $\theta$-functions.
The invariant to be computed was depending on an elliptic curve
(with modulus $q=e^{\pi \,i\,\tau}$) and several points on the curve. 
It appeared as a sum of $1440$ terms, each an integral over
a period of a rational fraction of high degree ($16$)
in $\theta$-functions and their derivatives. 
After computing the first terms in the $q$-expansion of the sum, 
(with the help of a computer\footnote{We wish to express our gratitude to
Michael Trott for his kind assistance}),
and factoring out basic elliptic functions of the above parameters,
we were left with a scalar function of $q$, starting as,\\
$$
q^{3/4} - 9\, q^{11/4} + 27\, q^{19/4} - 12\, q^{27/4} - 90\, q^{35/4} +\ldots 
$$
in which one recognises the $9^\mathrm{th}$ power of the Dedekind $\eta$-function.\\

\noindent We then gradually simplified the result (with $\eta^9$ appearing 
as an integration factor from the derivative of the Weierstrass $\wp$-function)
and elaborated the concepts which directly explain the final  form of the result.\\

\noindent The main new conceptual tool which we obtained and that we develop 
in this paper is the notion of central quadratic form
for quadratic algebras (Definition \ref{cent}). The 
 geometric data $\{E\,,\,  \sigma\,,\,\call\}$ of a quadratic algebra 
$\cala$ is a standard notion (\cite{ode-fei:1989}, \cite{art-tat-vdb:1990},
 \cite{smi-sta:1992}) defined in
such a way that the algebra maps homomorphically to a cross-product algebra
obtained from  sections of  powers of the line bundle $\call$ on 
powers of  the correspondence $\sigma$.
 We shall prove a purely algebraic result (Lemma \ref{alg0})
which considerably refines the above homomorphism and lands
in a richer cross-product. 
Its origin  can be traced as explained above
to the work of Sklyanin (\cite{skl:1982})
and  Odesskii-Feigin (\cite{ode-fei:1989}).\\

\noindent Our construction is then refined to control the $C^\ast$-norm 
(Theorem \ref{C*}) and  to construct the differential calculus in 
the desired generality (section 8).
It also allows  to show the pertinence of the general ``spectral"
framework for noncommutative geometry, in spite of the rather
esoteric nature of the above examples.\\

\noindent 
We  apply the general theory 
of central quadratic forms to show that the noncommutative 3-spheres $S^3_\ug$
 admit a natural ramified covering $\pi$ by a noncommutative
$3$-manifold $M$ which is (in the even case, cf. Proposition \ref{siginvar})
  isomorphic to the mapping torus of an outer automorphism
of the noncommutative
$2$-torus $T^2_\eta$. It is a noncommutative version of a nilmanifold
(Corollary \ref{H3}) with a natural action of the Heisenberg Lie algebra
$\frach_3$ and an invariant trace.
 This covering yields the differential calculus in its ``transcendental" form.

\noindent  Another important novel concept which plays a basic role in the 
present paper is the notion of Jacobian for a morphism of noncommutative spaces,
developed from basic ideas of noncommutative differential geometry \cite{ac:1982}
and expressed in terms of Hochschild homology (section 7). We 
 compute the Jacobian of the above ramified covering $\pi$
by pairing the direct image of the fundamental class of the
noncommutative $3$-dimensional nilmanifold $M$
with the Chern character of the defining unitary 
and obtain the answer as the product of a period (of
an  elliptic integral) by a rational function (Theorem
\ref{vol} and Corollary \ref{nct0}).
As explained above, we first computed these expressions
in terms of elliptic functions and modular
forms which led us in order to simplify the results to 
 extend the moduli space from the real to the
complex domain, and to formulate everything  in geometric terms. 

\noindent The leaves of the scaling foliation of the 
real moduli space then 
appear as the real parts of a net
of degree $4$ elliptic curves in $P_3(\mathbb C)$
having  $8$ points in common.
These elliptic curves  turn out to play a fundamental role and 
to be  closely related to the elliptic curves 
of the geometric data (cf. Section $4$) of the quadratic algebras which their elements label. 
In fact we first directly compared the $j$-invariant of the
generic fibers (of the scaling foliation) with the $j$-invariant of the 
quadratic algebras, and found them to be equal. 
This equality is surprising in that it fails in the degenerate
(non-generic) cases, where the characteristic variety can be as
large as $P_3(\mathbb C)$. Moreover, even in the generic
case, the two notions of ``real" points, either in the fiber
or in the characteristic variety are not the same, but dual
to each other. We eventually explain in Theorem
\ref{iden} the generic coincidence between 
the leaves of the scaling foliation in the complexified 
moduli space and the characteristic varieties of the associated
quadratic algebras. This Theorem
\ref{iden}  also exhibits the relation of our theory with iterations of
a specific birational automorphism of $P_3(\mathbb C)$, defined over $\mathbb Q$,
and restricting on each fiber as a translation of this elliptic curve.
The generic irrationality of this translation and the nature of 
its diophantine approximation play an important role in sections 7 and 8. \\

\noindent Another important result is the appearance of the hyperfinite 
factor of type II$_1$ as cross-product of the field $K_q$
of meromorphic functions on an elliptic curve by a subgroup of the Galois group 
Aut$_{\mathbb C}(K_q)$, (Theorem~\ref{gal}) and the description of the 
differential calculus in general (Lemma \ref{cyclic}) and 
in ``rational" form on $S^3_\ug$ 
(Theorem~\ref{fund}). 
The detailed proofs together with the 
analysis of the spectral geometry of $S^3_\ug$ 
and of the $C^\ast$-algebra $C^\ast(S^3_\ug)$ will appear in Part II.\\

\section{The Real Moduli Space of 3-Spheres $S^3_\ug$}

\noindent Let us now be more specific and describe the 
basic $K$-theoretic equation defining our spheres.
In the simplest
case it asserts that the algebra $\cala$ of ``coordinates" on the 
noncommutative space is generated by a self-adjoint idempotent $e$, ($e^2=e$, $e=e^{\ast}$)
together with the algebra $M_2(\mathbb C)$ of two by two scalar matrices.
 The only relation is
that the trace of $e$ vanishes, i.e. that the projection of $e$ on the commutant
of $M_2(\mathbb C)$ is zero. One shows that $\cala$ is then the algebra 
$M_2(\calg(S^2))$ where $\calg(S^2)$ is the algebra of coordinates on the standard
$2$-sphere $ S^2$.\\

\noindent The general form of the equation distinguishes two cases according to
the parity of the dimension $d$. In the even case 
the algebra $\cala$ of ``coordinates" on the 
noncommutative space is still generated by a projection $e$, ($e^2=e$, $e=e^{\ast}$)
and an algebra of scalar matrices, but the dimension $d = 2 k$ appears in
requiring the vanishing not only of the trace of $e$, but of all
components of its Chern character of degree $\,0, \ldots, d-2$.
This of course involves the algebraic (cyclic homology) formulation of
the Chern Character.

\noindent The Chern character in cyclic homology \cite{ac:1982}, \cite{ac:1986a}, 
\begin{equation}
\rm{ch}_* \,:\quad  K_*(\cala) \rightarrow HC_*(\cala)
\end{equation}
 is the noncommutative geometric analogue of the classical Chern character.
We describe it in the odd case which is relevant in our case $d=3$.
Given a noncommutative algebra $\cala$, an invertible element
 $U$ in $M_d(\cala)$ defines a class in $K_1(\cala)$ and the 
components of its Chern character are given by,
\begin{eqnarray}
\ch_{\frac{n}{2}}(U)&=& U^{i_0}_{i_1}\otimes V^{ i_{1}}_{i_2}
\otimes \cdot \cdot \cdot \otimes U^{i_{n-1}}_{i_n}\otimes V^{i_n}_{i_0}
\nonumber \\
 &-& V^{ i_0}_{i_1} \otimes U^{i_1}_{i_2}\otimes\cdot \cdot \cdot \otimes V^{ i_{n-1}}_{i_n}\otimes U^{i_n}_{i_0}
\end{eqnarray}
where $V:= U^{-1}$ and summation  over repeated indices is understood.

\noindent By a noncommutative $n$-dimensional spherical manifold (n odd),
we mean the noncommutative space $S$ dual to the $\ast$-algebra $\cala$ 
generated
by the components $U^{i}_{j}$ of a \underline{unitary} solution
$U \in M_d(\cala)$, $d=2^{\frac{n-1}{2}}$, of the equation
\begin{equation} \label{sph}
\ch_{\frac{k}{2}}(U)=0\, , \quad \forall k < n\,, \; k \,{\rm odd}\,, \quad 
\ch_{\frac{n}{2}}(U)\neq 0
\end{equation}
which is the noncommutative counterpart of the vanishing
of the lower Chern classes of the Bott generator of the 
$K$-theory of classical odd spheres.\\


 \noindent The moduli space of 3-dimensional spherical manifolds 
appears naturally as a quotient of the space of symmetric unitary matrices
\begin{equation}
\cals \,:= \{\, \Lambda \in M_4(\mathbb C)\,\vert \quad \Lambda=\Lambda^t
\, , \quad \Lambda^\ast=\Lambda^{-1}\, \}
\end{equation}
Indeed for any  $\Lambda \in  \cals$ let $U( \Lambda)$ 
be the unitary solution of (\ref{sph}) given by 
\begin{equation}
U=\bbbone_2 \otimes z^0+i\sigma_k \otimes z^k
\label{eq1.1}
\end{equation}
where $\sigma_k$ are the usual Pauli matrices and
the presentation of the involutive algebra $ C_{\mathrm{alg}}(S^3(\Lambda))$ 
generated by the $z^\mu$ is given by the relations
\begin{equation} \label{pres}
U^\ast\, U\, =\, U\,U^\ast\,=\,1\, , \quad z^{\mu\ast}=\Lambda^\mu_\nu z^\nu
\end{equation}
We let $\cala:= C_{\mathrm{alg}}(\mathbb R^4(\Lambda))$ be the associated
quadratic algebra, generated by the $z^\mu$ with presentation,
\begin{equation} \label{pres1}
U^\ast\, U\, =\, U\,U^\ast\,\in \bbbone_2 \otimes \cala \,, \quad z^{\mu\ast}=\Lambda^\mu_\nu z^\nu
\end{equation}
where $\bbbone_2$ is the unit of $M_2(\mathbb C)$.
The element $r^2 :=\, \sum^3_{\mu=0} z^{\mu} z^{\mu\ast}
=\, \sum^3_{\mu=0} z^{\mu\ast} z^\mu$ is 
in the center of $ C_{\mathrm{alg}}(\mathbb R^4(\Lambda))$ (cf. Part I) and 
the additional inhomogeneous relation defining $ C_{\mathrm{alg}}(S^3(\Lambda))$
is $r^2 =1$.

\noindent By Part I Theorem  1, any unitary solution 
of (\ref{sph}) for $n=3$
is a homomorphic image of $U( \Lambda)$ 
for some $ \Lambda \in \cals $.
\noindent The transformations
\begin{equation}
U\mapsto  \lambda \,U\, \, \,  \mbox{with}\, \, \, \lambda = e^{i\varphi} \in U(1)
\label{eq1.2}
\end{equation}
\begin{equation}
U\mapsto V_1 \, U \, V_2\, \,\,  \mbox{with}\, \,\,    V_1, V_2 \in SU(2)
\label{eq1.3}
\end{equation}
\begin{equation}
U\mapsto U^\ast
\label{eq1.4}
\end{equation}
act on the space of unitary solutions of (\ref{sph})
and preserve the isomorphism class of the algebra 
$ C_{\mathrm{alg}}(S^3(\Lambda))$   and of the associated quadratic
algebra $ C_{\mathrm{alg}}(\mathbb R^4(\Lambda))$.  Transformation (\ref{eq1.2})
corresponds to 
\begin{equation}
\Lambda \mapsto e^{-2i\varphi} \, \Lambda 
\label{scal}
\end{equation}
\noindent  Transformation (\ref{eq1.3}) induces 
 $z^\mu \mapsto S^\mu_\nu z^\nu$ with $S\in SO(4)$ which in turn induces 
\begin{equation}
\Lambda \mapsto S\,\Lambda\, S^t \,, \quad S\in SO(4)
\label{rot}
\end{equation}
Finally (\ref{eq1.4}) reverses the ``orientation" of 
$S^3(\Lambda)$ and corresponds to 
\begin{equation}
\Lambda \mapsto \Lambda^{-1} 
\label{inv}
\end{equation}
We define the {\it real moduli space } $\calm$ as the quotient 
 of $\cals $ by the transformations (\ref{scal}) and 
(\ref{rot}), and the {\it unoriented} real moduli space $\calm^{\,'}$
as its quotient by (\ref{inv}).

\noindent By construction the space $\cals$ is
 the homogeneous space $U(4)/O(4)$, with $U(4)$ acting on $\cals$ by
\begin{equation}
\Lambda \mapsto V\Lambda V^t
\label{eq2.1}
\end{equation}
for $\Lambda \in \cals$ and $V\in U(4)$.  The conceptual description of 
the moduli spaces $\calm$ and $\calm^{\,'}$
requires taking care of finer details. We let $\theta$ be the involution of 
the compact Lie group $SU(4)$ given by complex conjugation,
and define the closed subgroup $K \subset SU(4)$ as the normaliser 
of $ SO(4) \subset SU(4)$, i.e.
by the condition 
\begin{equation}
K := \, \{u \in SU(4)\,\vert \, u^{-1}\, \theta(u) \, \in\, Z \}
\label{kk}
\end{equation}
where $Z$ is the center of $SU(4)$. The 
quotient $X:= SU(4)/K$ is a Riemannian globally symmetric space
(cf. \cite{hel:1978} Theorem 9.1, Chapter VII). One has $Z \subset K$ 
but the image of $K$ in 
$U:=SU(4)/Z$ is disconnected. 
Indeed, besides $Z \cdot SO(4)$ the subgroup $K$ contains the 
diagonal matrices with $\{v,v,v,v^{-3}\}$ as diagonal
elements, where $v$ is an $8$th root of $1$.

 \begin{proposition} \label{sym1}
The real moduli space $\calm$ (resp. $\calm^{\,'}$)
is canonically isomorphic to the space 
of congruence classes of point pairs
under the action of $SU(4)$
(resp. of isometries) in the Riemannian globally
symmetric space $X= SU(4)/K$.
\end{proposition}

\noindent This gives two equivalent descriptions of $\calm$
as the orbifold quotient of a $3$-torus by the action
of the Weyl group of the symmetric pair. 

\noindent In the first 
$(A_3)$ we identify the Lie algebra $\fracsu(4)=\Lie (SU(4))$
   with the Lie algebra of traceless antihermitian elements of $M_4(\mathbb C)$. We let $\fracd$ be
the Lie subalgebra of diagonal matrices, it is a maximal abelian subspace of $\fracsu(4)_-=\{ X\in \fracsu (4)\vert \theta(X)=-X\}$. The Weyl group $W$ of the symmetric pair $(SU(4),K)$ is isomorphic to the permutation group $\fracS_4$
acting on $\fracd$ by permutation of the matrix elements. 
The $3$-torus $T_A$ is the quotient,
\begin{equation} \label{ta}
T_A \, := \, \fracd \, / \Gamma\, , \quad 
\Gamma=\{ \delta \in \fracd\vert  e^{\delta}  \in K \}
\end{equation}
and the isomorphism of $T_A \,/ W$ with $\calm$ is obtained from, 
\begin{equation} \label{ta1}
\sigma(\delta)\, :=\,  e^{2 \delta}  \in \cals \, 
\quad \forall \delta \in \fracd
\end{equation}
The lattice $\Gamma$ is best expressed as $\Gamma = \{ \delta \in \fracd \vert \langle \delta,\rho\rangle \in \pi i \,\mathbb Z,\, \, \forall \rho \in \Delta\}$ in terms of the roots $\rho$ of the pair $(SU(4),K)$ where the root system $\Delta$ is the same as for the Cartan subalgebra $\fracd_\mathbb C=\fracd \otimes \mathbb C \subset \fracsl(4,\mathbb C)$ of $\fracsl (4,\mathbb C)=\fracsu(4)_\mathbb C$. The roots $\alpha_{\mu,\nu}$, $\mu,\nu\in \{0,1,2,3\}$, $\mu\not=\nu$ are given by
\begin{equation}
\alpha_{\mu,\nu}(\delta) = \alpha_\mu(\delta)-\alpha_\mu(\delta)
\label{eq3.9}
\end{equation}
where $\alpha_\mu(\delta)$ for $\mu\in \{0,1,2,3\}$ are the elements of the diagonal matrix $\delta\in \fracd$. In terms of the primitive roots
$\alpha_{0,k},\, \,  k\in \{1,2,3\}$
the coordinates $\varphi_k$ used in Part I are given  by $\varphi_k=\frac{1}{i} \alpha_{0,k}$
and the lattice $\Gamma$ corresponds to $\varphi_k\in \pi\mathbb Z$, $\forall k \in \{1,2,3\}$. With $T:= \mathbb R /\pi\mathbb Z$ we let  
 $\delta : T^3\, \mapsto  T_A $ be the inverse isomorphism.
 One has, modulo projective equivalence,
 \begin{equation} \label{phi-a}
\sigma(\delta_\varphi)\,\sim
\left[
\begin{array}{cccc}
1& & &\\
& e^{-2i\varphi_1} & &\\
&& e^{-2i\varphi_2}&\\
& & & e^{-2i\varphi_3}
\end{array}
\right]
\end{equation}
In these coordinates ($\varphi_k$) the symmetry given by the Weyl group $W=\fracS_4$ of the symmetric pair $(SU(4),K)$ now reads as follows,
\begin{equation}
\begin{array}{lll}
T_{01}(\varphi_1,\varphi_2, \varphi_3) & = & (-\varphi_1,\varphi_2-\varphi_1, \varphi_3-\varphi_1)\\
\\
T_{12}(\varphi_1,\varphi_2, \varphi_3) & = & (\varphi_2,\varphi_1, \varphi_3)\\
\\
T_{23}(\varphi_1,\varphi_2, \varphi_3) & = & (\varphi_1,\varphi_3, \varphi_2)
\end{array}
\label{eq2.2}
\end{equation}
where $T_{\mu\nu}$ is the transposition of $\mu,\nu\in \{0,1,2,3\}$, $\mu<\nu$.

\noindent In Part I we used 
the parametrization by $\ug\in T^3$ to label the 
 3-dimensional spherical manifolds $S^3_\ug$
and their quadratic counterparts $\mathbb R^4_\ug$.
The presentation of the algebra $C_{\mathrm{alg}}(\mathbb R^4_\ug)$
is given as follows, using the selfadjoint generators 
$x^\mu=x^{\mu\ast}$ for $\mu\in \{0,1,2,3\}$ related to the $z^k$
by $z^0=x^0$, $z^k=e^{i\varphi_k}$ $x^k$ for $j\in \{1,2,3\}$
\begin{equation}
\cos (\varphi_k) [x^0,x^k]_-=i \,\sin (\varphi_\ell-\varphi_m) [x^\ell, x^m]_+
\label{eq1.6}
\end{equation}
\begin{equation}
\cos(\varphi_\ell-\varphi_m)[x^\ell, x^m]_-=-i \,\sin (\varphi_k)[x^0,x^k]_+
\label{eq1.7}
\end{equation}
for $k=1,2,3$ where ($k, \ell, m$) is the cyclic permutation of ($1,2,3$) starting with $k$ and where $[a,b]_\pm=ab\pm ba$. 
The algebra $C_{\mathrm{alg}}( S^3_\ug)$ is the quotient of 
$C_{\mathrm{alg}}(\mathbb R^4_\ug)$  by the two-sided ideal generated by the hermitian central element $\sum_\mu(x^\mu)^2-\bbbone$. 

\noindent The second equivalent description of $\calm$
relies on the equality $A_3=D_3$, i.e. the identification of 
$SU(4)$ with the Spin covering of $SO(6)$
using the ($4$-dimensional) half spin representation of the latter.
Proposition \ref{sym1} holds unchanged replacing the pair
($SU(4), K$) by the pair ($SO(6), K_D$), where $K_D$
corresponds under the above isomorphism with the quotient
of $K$ by the kernel of the covering Spin$(6) \mapsto SO(6)$.
Identifying the Lie algebra $\fracso(6)$ with the Lie algebra
of antisymmetric six by six real matrices, $\fracd$
corresponds to the subalgebra $\fracd_D$ of
block diagonal matrices with three
blocks of the form 
\begin{equation}
\begin{pmatrix} 0 &\psi_k \\
 -\psi_k &0 \end{pmatrix}\end{equation}
The $3$-torus $T_D$ is the quotient,
\begin{equation} \label{td}
T_D \, := \, \fracd_D \, / \Gamma_D\, , \quad 
\Gamma_D=\{ \delta \in \fracd_D \,\vert  e^{\delta}  \in K_D \}
\end{equation}
which using the roots of $D_3$, $(\pm \, e_i \; \pm \, e_j)(\psi) :=\pm \, \psi_i \; \pm \, \psi_j
$ becomes, 
\begin{equation}
 \Gamma_D=\{ \psi\,  \vert\, \psi_i \; \pm \, \psi_j \in \pi\mathbb Z \}
\end{equation}

\noindent The relation between the parameters ($\psi_k$) for $k\in \{1,2,3\}$ and
the ($\varphi_k$) is given by,
\begin{equation}
2\,\psi_1 = \varphi_2+\varphi_3-\varphi_1,\,\,2\, \psi_2 = \varphi_3+\varphi_1-\varphi_2,\,\, 2\, \psi_3 = \varphi_1+\varphi_2-\varphi_3
\label{eq2.5}
\end{equation}
and in terms of the $\psi_k$ the Killing metric on $T_D$ reads,
\begin{equation}
ds^2=(d\psi_1)^2+(d\psi_2)^2+(d\psi_3)^2
\label{eq2.10}
\end{equation}
The action of the Weyl group $W$ is given by the subgroup of $O(3,\mathbb Z)$
\begin{equation}
w(\psi_1,\psi_2,\psi_3)=(\varepsilon_1\psi_{\sigma(1)},\varepsilon_2 \psi_{\sigma(2)},\varepsilon_3\psi_{\sigma(3)})
\label{eq2.7}
\end{equation}
with $\sigma\in \fracS_3$ and $\varepsilon_k\in \{1,-1\}$,  of  elements  $w \in O(3,\mathbb Z)$
such that $\varepsilon_1\varepsilon_2\varepsilon_3=1$.
The additional symmetry (\ref{inv})
defining $\calm^{\,'}$ is simply 
\begin{equation}
  \psi\,  \mapsto -\psi
\end{equation}
and together with $W$ it generates $ O(3,\mathbb Z)$.

 \noindent Another advantage of the variable $\psi_k$ is that the scaling vector field $Z$ (cf. Part I) whose orbits describe the local equivalence 
relation on $\calm$ generated by the isomorphism of quadratic algebras,
\begin{equation}
C_{\mathrm{alg}}(\mathbb R^4(\Lambda_1)) \sim C_{\mathrm{alg}}(\mathbb R^4(\Lambda_2))
\end{equation}
and which was given in the variables ($\varphi_k$) as,
\begin{equation}
Z=\sum^3_{k=1} \sin (2\varphi_k)\sin (\varphi_\ell + \varphi_m -\varphi_k)\frac{\partial}{\partial \varphi_k}
\label{eq1.9}
\end{equation}
  is now given by
\begin{equation}
Z=\frac{1}{4}\, \sum^3_{k=1} \frac{\partial H_0}{\partial \psi_k} \frac{\partial}{\partial \psi_k}
\label{eq2.8}
\end{equation}
with
\begin{equation}
H_0=\sin(2\,\psi_1)\sin(2\,\psi_2)\sin(2\,\psi_3)\,.
\label{eq2.9}
\end{equation}

\noindent Since $2\, \Gamma_D$ is the unit lattice for $PSO(6)$
we can translate everything to the space $\cal{C}$ of conjugacy classes 
in $PSO(6)$.\\

\begin{theorem}\label{scaltheo}
The doubling map $\psi \mapsto 2 \psi$ establishes an 
isomorphism between $\calm$ and $\cal{C}$ transforming the scaling foliation on $\calm$ into the gradient flow
(for the Killing metric) of the character of the signature representation of $SO(6)$, i.e. the super-trace of its action on 
$\wedge^3\mathbb C^6 =\wedge^3_+\mathbb C^6 \oplus \wedge^3_-\mathbb C^6 $ with    $\wedge^3_\pm\mathbb C^6   = \{\omega\in \wedge^3\mathbb C^6 \vert \ast \omega=\pm \, i\omega\}$. 
\end{theorem}

\noindent This  flow admits remarkable compatibility properties
with the canonical cell decomposition of $\cal{C}$, they will be analysed in
Part II.

\bigskip

\section{The Complex Moduli Space and its Net of Elliptic Curves}

\noindent The proper understanding of the noncommutative spheres $S_\ug^3$ 
relies on basic computations (Theorem \ref{vol} below) whose result
depends on $\ug$ through elliptic integrals.
The conceptual explanation of this dependence requires  
 extending the moduli space from the real to the
complex domain. 
The leaves of the scaling foliation then 
appear as the real parts of a net
of degree $4$ elliptic curves in $P_3(\mathbb C)$
having  $8$ points in common.
These elliptic curves will turn out to play a fundamental role and 
to be  closely related to the elliptic curves 
of the geometric data (cf. Section $4$) of the quadratic algebras which their elements label. 

\noindent To extend the moduli space to the 
complex domain we start with the relations 
defining  the involutive algebra $ C_{\mathrm{alg}}(S^3(\Lambda))$ 
and take for $\Lambda$ the diagonal
matrix with 
\begin{equation}
\Lambda^\mu_\mu \, := \, u_\mu^{-1}
\end{equation}
where $(u_0,u_1,u_2,u_3)$ are the coordinates of
 $\ug\in  P_3(\mathbb C)$.
Using $y_{\mu}:=\Lambda^\mu_\nu z^\nu $ one obtains the 
homogeneous defining relations in the form,
\begin{eqnarray}
u_k\, y_k\, y_0-u_0\, y_0\, y_k+u_\ell\,  y_\ell\,  y_m-u_m\, y_m\, y_\ell & = & 0   \nonumber\\
u_k\, y_0\, y_k-u_0\, y_k\, y_0+u_m\,  y_\ell\,  y_m-u_\ell\,  y_m\, y_\ell & = & 0\label{homc}
\end{eqnarray}
for any cyclic permutation $(k,\ell,m)$ of (1,2,3).
The inhomogeneous relation becomes,
\begin{equation}
\sum \, u_\mu\, y_\mu^2 =\, 1
\end{equation}
and the corresponding algebra $ C_{\mathrm{alg}}(S_{\mathbb C}^3(\ug))$
only depends upon the class of $\ug\in  P_3(\mathbb C)$.
We let $ C_{\mathrm{alg}}(\mathbb C^4(\ug))$ be the 
quadratic algebra defined by (\ref{homc}). 

\noindent 
Taking $u_\mu\,= e^{2 i\,\varphi_\mu}$, $\varphi_0=0$,
 for all $\mu$ and $x_\mu:=e^{ i\,\varphi_\mu}y_\mu$
we obtain the defining relations
(\ref{eq1.6}) and (\ref{eq1.7}) 
 (except for $x_\mu^{\ast}=x_\mu$).

\noindent 
 We showed in Part I that 
for $\varphi_k\not= 0$ and $\vert\varphi_r-\varphi_s \vert\not=\frac{\pi}{2}$ (i.e. $u_k\not=1$ and $u_r\not= -u_s$) for any $k,r,s\in \{1,2,3\}$, one can find 4 scalars $s^\mu$ such that by setting $S_\mu=s^\mu x^\mu$ for $\mu\in \{0,1,2,3\}$ the system (\ref{eq1.6}), (\ref{eq1.7}) reads
\begin{equation}
[S_0,S_k]_-=iJ_{\ell m} [S_\ell, S_m]_{+}
\label{eq1.10}
\end{equation}
\begin{equation}
[S_\ell,S_m]_-=i[S_0,S_k]_{+}
\label{eq1.11}
\end{equation}
for $k=1,2,3$, where $(k,\ell,m)$ is the cyclic permutation of (1,2,3) starting with $k$ and where $J_{\ell m}=-\tan (\varphi_\ell-\varphi_m) \tan (\varphi_k)$. One has
\begin{equation}
J_{12} + J_{23} + J_{31} + J_{12}J_{23}J_{31}=0
\label{eq1.12}
\end{equation}
and the relations (\ref{eq1.10}), (\ref{eq1.11}) together with (\ref{eq1.12}) for the scalars $J_{\ell m}$ characterize the Sklyanin algebra \cite{skl:1982}, \cite{skl:1983}, a regular algebra of global dimension 4 which has been widely studied (see e.g. \cite{ode-fei:1989}, \cite{smi-sta:1992}, \cite{tat-vdb:1996}) and which plays an important role in noncommutative algebraic geometry. From (\ref{eq1.12}) it follows that for the above (generic) values of $\ug$, $C_{\mathrm{alg}}(\mathbb R^4_\ug)$
 only depends on 2 parameters;
with $Z$ as in (\ref{eq1.9})
 one has $Z(J_{k\ell})=0$ and the leaf of the scaling foliation through a generic $\ug \in T^3$ is the
connected component of $\ug$ in
 \begin{equation} \label{def0}
F_{\mathbb T}(\ug):= \{\vg\in T^3\ \vert\, J_{k\ell}(\vg)=\,
 J_{k\ell}(\ug)\}
\end{equation}

\noindent In terms of homogeneous parameters 
the functions $J_{\ell m}$ read as 
\begin{equation}
J_{\ell m}=\tan (\varphi_0-\varphi_k)\tan (\varphi_\ell-\varphi_m)
\end{equation}
for any cyclic permutation $(k,\ell, m)$ of (1,2,3),
and extend to the complex domain $\ug\in P_3(\mathbb C)$ as, 
\begin{equation}\label{jlm}
J_{\ell m}=\frac{(u_0+u_\ell)(u_m+u_k)-(u_0+u_m)(u_k+u_\ell)}{(u_0+u_k)(u_\ell+u_m)}
\end{equation}
It follows easily from the finer Theorem \ref{iden} that for generic
values of $\ug\in P_3(\mathbb C)$
the quadratic algebra $ C_{\mathrm{alg}}(\mathbb C^4(\ug))$
only depends upon $J_{k\ell}(\ug)$. We thus define
\begin{equation}  
F(\ug):= \{\vg\in P_3(\mathbb C)\, \vert\, J_{k\ell}(\vg)=\,
 J_{k\ell}(\ug)\}  \label{def}
\end{equation}
Let then, 
\begin{equation}
(\alpha,\beta, \gamma)=\left\{(u_0+u_1)(u_2+u_3),(u_0+u_2)(u_3+u_1),(u_0+u_3)(u_1+u_2)\right\}
\nonumber
\end{equation}
be the Lagrange resolvent of the 4th degree equation, 
\begin{equation}
\Phi(\ug)=(\alpha,\beta,\gamma)
\label{eq6.2}
\end{equation}
viewed as a map
\begin{equation}
\Phi:P_3(\mathbb C)\backslash S \rightarrow P_2(\mathbb C)
\label{eq6.5}
\end{equation}
where $S$ is the following set 
 of 8 points
\begin{eqnarray}
p_0=(1,0,0,0),\, \, p_1=(0,1,0,0),\, \, p_2=(0,0,1,0),\, \, p_3=(0,0,0,1)\quad \quad
\label{eq6.3}\\
q_0=(-1,1,1,1),\, \, q_1=(1,-1,1,1),\, \, q_2=(1,1,-1,1),\, \, q_3=(1,1,1,-1)
\nonumber
\end{eqnarray}

\noindent We extend the generic definition (\ref{def}) to arbitray
 $\ug\in P_3(\mathbb C)\backslash S$ and define  $F_\ug$ 
in general as
the union of $S$ with the fiber of $\Phi$ through $\ug$. It can be understood geometrically as follows.

\noindent Let $\caln$ be the net of quadrics in $P_3(\mathbb C)$ which contain $S$. Given $\ug\in P_3(\mathbb C)\backslash S$ the elements of $\caln$ which contain $\ug$ form a pencil of quadrics with base locus
\begin{equation}
\cap \{Q\, \vert \, Q\in \caln,\, \ug\in Q\}=Y_{\ug}
\label{eq6.6}
\end{equation}
which is an elliptic curve of degree 4 containing $S$ and $\ug$. One has
\begin{equation}
Y_{\ug}
=F_{\ug}
\label{eq6.7}
\end{equation}
$$  
\hbox{
\psfig{figure=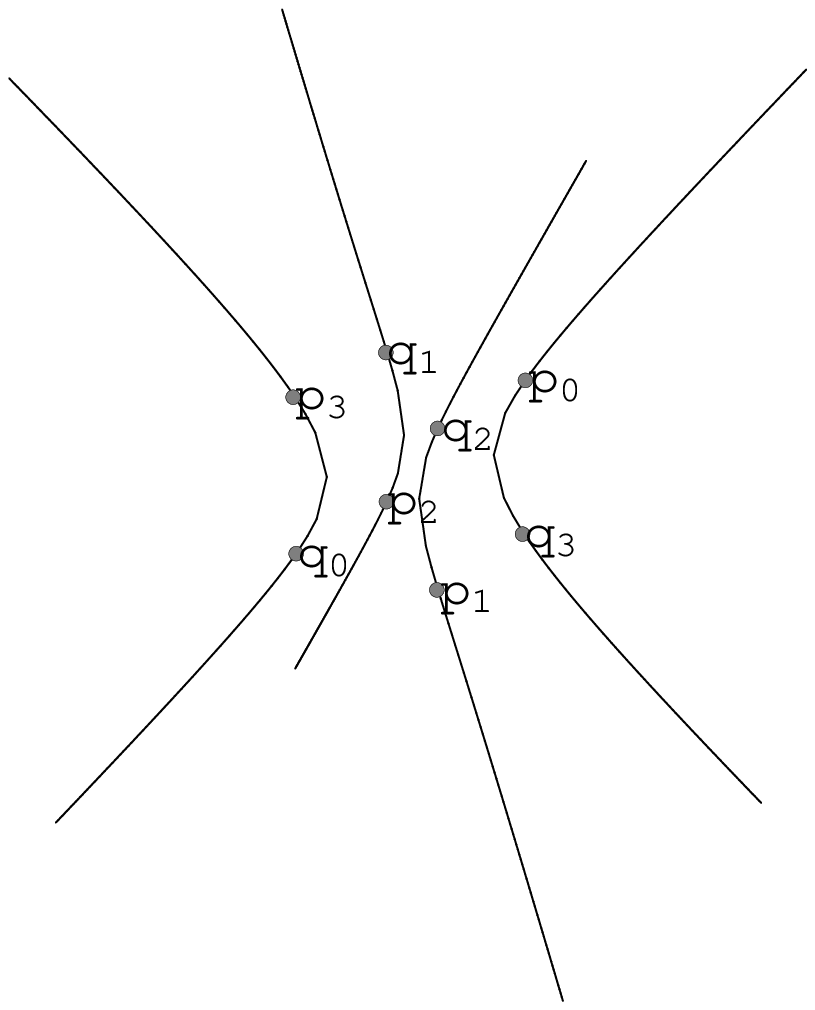}
}
$$
\centerline{Figure $1$: The Elliptic Curve $F_\ug \cap P_3(\mathbb {R})$ }

\noindent We shall now give, for generic values of $(\alpha,\beta,\gamma)$ a parametrization of $F_\ug$ by $\theta$-functions. We start with the equations for $F_\ug$
\begin{equation}
\frac{(u_0+u_1)(u_2+u_3)}{\alpha}=\frac{(u_0+u_2)(u_3+u_1)}{\beta}=\frac{(u_0+u_3)(u_1+u_2)}{\gamma}
\label{eq6.8}
\end{equation}
and we diagonalize the above quadratic forms as follows
\begin{equation}
\begin{array}{lll}
(u_0+u_1)\,(u_2+u_3)& = & Z_0^2-Z_1^2\\ 
(u_0+u_2)\,(u_3+u_1) & = & Z_0^2-Z_2^2\\
 (u_0+u_3)\,(u_1+u_2)& = &Z_0^2-Z_3^2
\end{array}
\label{eq6.9}
\end{equation}
where 
\begin{equation}
(Z_0,Z_1,Z_2,Z_3) = M.u 
\label{eq6.10}
\end{equation}
where $M$ is the involution,
\begin{equation}
M:= \frac{1}{2}\left[
\begin{array}{cccc}
1& \;1& 1&1\\
1& 1 &-1 &-1\\
1&-1& 1&-1\\
1& -1& -1& 1
\end{array}
\right]
\end{equation}
In these terms the equations for $F_\ug$ read
\begin{equation}
\frac{Z_0^2-Z_1^2}{\alpha}=\frac{Z_0^2-Z_2^2}{\beta}=\frac{Z_0^2-Z_3^2}{\gamma}
\label{eq6.11}
\end{equation}
Let now $\omega\in \mathbb C$, $\im \,\omega >0$ and $\eta\in \mathbb C$ be such that one has, modulo projective equivalence,
\begin{equation}
(\alpha, \beta, \gamma)\sim \left(\frac{\theta_2(0)^2}{\theta_2(\eta)^2},\, \, \frac{\theta_3(0)^2}{\theta_3(\eta)^2},\, \, \frac{\theta_4(0)^2}{\theta_4(\eta)^2}\right)
\label{eq6.12}
\end{equation}
where $\theta_1,\theta_2, \theta_3, \theta_4$ are the theta functions for the lattice $L=\mathbb Z +\mathbb Z\omega\subset \mathbb C$
\begin{proposition} \label{theta}
The following define isomorphisms of $\mathbb C/L$ with $F_\ug$,
\[
\varphi(z)=\left( \frac{\theta_1(2z)}{\theta_1(\eta)},\, \frac{\theta_2(2z)}{\theta_2(\eta)},\, \frac{\theta_3(2z)}{\theta_3(\eta)},\, \frac{\theta_4(2z)}{\theta_4(\eta)}\right) = (Z_0,Z_1,Z_2,Z_3)
\]
and $\psi(z)=\varphi(z-\eta/2)$. 
\end{proposition}
\noindent {\bf Proof}\hspace{0,2cm} Up to an affine transformation, $\varphi$ (and $\psi$ are) is the classical projective embedding of $\mathbb C/L$ in $P_3(\mathbb C)$. Thus we only need to check that the biquadratic curve $\im\,\varphi=\im\,\psi$ is given by (\ref{eq6.11}). It is thus enough to check (\ref{eq6.11}) on $\varphi(z)$. This follows from the basic relations
\begin{equation}
\theta^2_2(0)\theta^2_3(z)=\theta^2_2(z)\theta^2_3(0)+\theta^2_4(0)\theta^2_1(z)
\label{eq6.13}
\end{equation}
and
\begin{equation}
\theta^2_4(z)\theta^2_3(0)=\theta^2_1(z)\theta^2_2(0)+\theta^2_3(z)\theta^2_4(0)
\label{eq6.14}
\end{equation}
which one uses to check $\frac{Z_0^2-Z_1^2}{\alpha}=\frac{Z_0^2-Z_2^2}{\beta}$ and $\frac{Z_0^2-Z_2^2}{\beta}=\frac{Z_0^2-Z_3^2}{\gamma}$ respectively.$\square$\\

\noindent The elements of $S$ are obtained from the following values of $z$
\begin{equation}
\psi(\eta)=p_0,\, \psi(\eta+\frac{1}{2})=p_1,\, \psi(\eta+\frac{1}{2}+\frac{\omega}{2})=p_2,\, \psi(\eta+\frac{\omega}{2})=p_3
\label{eq6.15}
\end{equation}
and
\begin{equation}
\psi(0)=q_0,\, \psi(\frac{1}{2})=q_1,\, \psi(\frac{1}{2}+\frac{\omega}{2})=q_2,\, \psi(\frac{\omega}{2})=q_3.
\label{eq6.16}
\end{equation}
(We used $M^{-1}.\psi$ to go back to the coordinates $u_\mu$).\\
Let $H\sim \mathbb Z_2\times \mathbb Z_2$ be the Klein subgroup of the symmetric group $\fracS_4$ acting on $P_3(\mathbb C)$ by permutation of the coordinates $(u_0,u_1,u_2,u_3)$.\\
For $\rho$ in $H$ one has $\Phi\circ \rho=\Phi$, so that $\rho$ defines for each $\ug$ an automorphism of $F_\ug$. For $\rho$ in $H$ the matrix $M\rho M^{-1}$ is diagonal with $\pm 1$ on the diagonal and the quasiperiodicity of the $\theta$-functions  allows to check that these automorphisms are translations on $F_\ug$ by the following 2-torsion elements of $\mathbb C/L$,
\begin{equation}
\begin{array}{ll}
\rho =\left(
\begin{array}{cccc}
0 & 1 & 2 & 3\\
1 & 0 & 3 & 2
\end{array}\right) & \mbox{is translation by}\, \frac{1}{2},\\
\\
\rho =\left(
\begin{array}{cccc}
0 & 1 & 2 & 3\\
2 & 3 & 0 & 1
\end{array}\right) & \mbox{is translation by}\, \frac{1}{2}+\frac{\omega}{2}
\end{array}
\label{eq6.17}
\end{equation}
Let $\calo\subset P_3(\mathbb C)$ be the complement of the $4$ hyperplanes $\{u_\mu\,=0\}$ with $ \mu\in \{0,1,2,3\}$. Then $(u_0,u_1,u_2,u_3)\mapsto (u_0^{-1},u_1^{-1},u_2^{-1},u_3^{-1})$ defines an involutive automorphism $I$ of $\calo$ and since one has
\begin{equation}
(u_0^{-1}+u_k^{-1})(u_\ell^{-1}+u^{-1}_m)=(u_0u_1u_2u_3)^{-1}(u_0+u_k)(u_\ell+u_m)
\label{eq6.18}
\end{equation}
it follows that $\Phi\circ I=\Phi$, so that $I$ defines for each $\ug\in \calo\backslash \{q_0,q_1,q_2,q_3\}$ an involutive automorphism of $F_\ug\cap\calo$ which extends canonically to $F_\ug$, in fact,
\begin{proposition}\label{nice}
The restriction of $I$ to $F_\ug$ is the symmetry
$\psi(z)\mapsto \psi(-z)$ around any of the
points $q_\mu \in F_\ug$
in the elliptic curve $F_\ug$.
\end{proposition}
This symmetry, as well as the above translations 
by two torsion elements does not refer to a choice of origin
in the curve $F_\ug$.
The proof follows from identities on theta functions.

\noindent Let $T  :=\{\ug\,\vert \, \vert u_\mu \vert
 =1  \; \forall \mu \}$. By section $2$ the torus
$T$ gives a covering of the real moduli space $\calm$. For $\ug \in T$, the 
point  $\Phi(\ug)$ is real with projective coordinates
\begin{equation}\label{real0}
\Phi(\ug) = (s_1,s_2,s_3)\, , \quad s_k := 1 + t_\ell\, t_m
\, , \quad t_k:= {\rm tan}(\varphi_k- \varphi_0)
\end{equation}
The corresponding fiber $F_{\ug}$ is stable under 
complex conjugation $\vg \mapsto \overline{\vg}$
and the intersection of $F_{\ug}$ with the real moduli space is 
given by,
\begin{equation}\label{real}
F_{\mathbb T}(\ug)=\, F_{\ug}\cap T=\{\vg \in F_{\ug} \vert I(\vg)\,=\, \overline{\vg}\}
\end{equation}
The curve $F_{\ug}$ is defined over $\mathbb R$ and (\ref{real})
determines its purely \underline{imaginary} points. 
 Note that  $F_{\mathbb T}(\ug) $ (\ref{real})
is  invariant under the Klein group $H$ 
and thus has two connected components, we let 
$F_{\mathbb T}(\ug)^0 $ be the component containing $q_0$.
The real points,
$ \{\vg \in F_{\ug} \vert \vg\,=\, \overline{\vg}\}=\,F_\ug \cap P_3(\mathbb {R})
$ of $F_{\ug}$ do play a complementary role in the characteristic variety
(Proposition \ref{siginvar}).

\section{ Generic Fiber $=$ Characteristic Variety}

\noindent Let us recall the definition of the geometric data $\{E\,,\,  \sigma\,,\,\call\}$ for  quadratic  algebras. Let $\cala=A(V,R)=T(V)/(R)$ be a quadratic algebra where $V$ is a finite-dimensional complex vector space and where $(R)$ is the two-sided ideal of the tensor algebra $T(V)$ of $V$ generated by the subspace $R$ of $V\otimes V$. Consider the subset of $V^\ast \times V^\ast$ of  pairs $(\alpha,\beta)$ such that
\begin{equation}
\langle \omega,\alpha\otimes\beta\rangle=0,\, \, \, \alpha\not= 0, \beta\not=0
\label{eq5.1}
\end{equation}
for any $\omega\in R$. Since $R$ is homogeneous, (\ref{eq5.1}) defines a subset 
\[
\Gamma \subset P(V^\ast)\times P(V^\ast)
\]
 where $P(V^\ast)$ is the complex projective space of one-dimensional complex subspaces of $V^\ast$. Let $E_1$ and $E_2$ be the first and the second projection of $\Gamma$ in $P(V^\ast)$. It is usually assumed that they coincide i.e. that one has
\begin{equation}
E_1=E_2=E\subset P (V^\ast)
\label{eq5.2}
\end{equation}
and that the correspondence $\sigma$ with graph $\Gamma$ is an automorphism of $E$, $\call$ being the pull-back on $E$ of the dual of the tautological line bundle of $P(V^\ast)$. The algebraic variety $E$ is refered to as the characteristic variety. In many cases $E$ is the union of an elliptic curve with a finite number of points which are invariant by $\sigma$. This is the case for $\cala_\ug= C_{\mathrm{alg}}(\mathbb C^4(\ug))$  at generic $\ug$ since  it then reduces to the Sklyanin algebra for which this is known \cite{ode-fei:1989}, \cite{smi-sta:1992}. When $\ug$ is non generic, e.g. when the isomorphism with the Sklyanin algebra breaks down, the situation is more involved, and the characteristic variety can be as large as $P_3(\mathbb C)$. This is described in Part I, where we gave a complete description of the geometric datas. Our aim now is to show that for $\ug\in P_3(\mathbb C)$ generic, there is an astute  choice of generators of the quadratic algebra $\cala_\ug= C_{\mathrm{alg}}(\mathbb C^4(\ug))$  for which  the characteristic variety $E_\ug$ actually coincides with the fiber variety $F_\ug$ and to identify the corresponding automorphism $\sigma$. Since this coincidence  no longer holds for non-generic values it is a quite remarkable fact which we first noticed by comparing the $j$-invariants of these two elliptic curves.\\

\noindent Let  $\ug\in P_3(\mathbb C)$ be generic, 
we perform the following change of generators

\begin{equation} \label{change0}
\begin{array}{llll} 
y_0 & = & \sqrt{u_1-u_2} \; \sqrt{u_2-u_3}  \;\sqrt{u_3-u_1} & Y_0\\
\\
y_1 & = & \sqrt{u_0+u_2} \; \sqrt{u_2-u_3}  \;\sqrt{u_0+u_3}&  Y_1\\
\\
y_2 & = & \sqrt{u_0+u_3}  \;\sqrt{u_3-u_1}  \;\sqrt{u_0+u_1}&  Y_2\\
\\
y_3 & = & \sqrt{u_0+u_1}  \;\sqrt{u_1-u_2}  \;\sqrt{u_0+u_2}&  Y_3
\end{array}
\end{equation}

\noindent We let $J_{\ell m}$ be as before, given by (\ref{jlm})
\begin{equation}
J_{12}=\frac{\alpha-\beta}{\gamma},\, \, \, J_{23}=\frac{\beta-\gamma}{\alpha},\, \, \,  J_{31}=\frac{\gamma -\alpha}{\beta}
\label{eq5.8}
\end{equation}
with $\alpha, \beta, \gamma$ given by (\ref{eq6.2}).
Finally let $e_\nu$ be the $4$ points of $P_3(\mathbb C)$
whose homogeneous coordinates ($Z_\mu$) all vanish but one. \\

\begin{theorem} \label{iden}
\noindent 1) In terms of the $Y_\mu$, the relations of $\cala_\ug$ take the form 
\begin{eqnarray}
[Y_0,Y_k]_- & = & [Y_\ell,Y_m]_+\label{eq5.6}\\
{[Y_\ell,Y_m]}_- & = & - J_{\ell m} [Y_0,Y_k]_+\label{eq5.7}
\end{eqnarray}
for any $k\in \{1,2,3\}$, $(k,\ell,m)$ being a cyclic permutation of (1,2,3)

\noindent 2) The characteristic variety $E_\ug$ is the union of $F_\ug$ with the $4$ points
$e_\nu$.

\noindent 3) The automorphism $\sigma$ of the characteristic variety $E_\ug$  is given by 
\begin{equation}
\psi(z) \mapsto \psi(z- \eta)
\label{eq5.10}
\end{equation}
on $F_\ug$ and $\sigma=\id$ on the 4 points $e_\nu$.

\noindent 4) The automorphism $\sigma$ is the restriction to $F_\ug$ of a birational
automorphism of $P_3(\mathbb C)$ independent of $\ug$ and defined over $\mathbb Q$.
\end{theorem} 

\noindent The resemblance between the above presentation and the Sklyanin one (\ref{eq1.10}), (\ref{eq1.11}) is misleading, for the latter all the 
characteristic varieties are contained in the same quadric (cf. \cite{smi-sta:1992} \S 2.4)
$$
\sum x_\mu^2= 0
$$
and cant of course form a net of essentially disjoint curves.

\noindent {\bf Proof} $\;\;$ By construction $E_\ug=\{ Z \, \vert \, {\rm Rank}\;N(Z)<4 \}$ where

\medskip
\begin{equation}
N(Z)=\left(
\begin{array}{cccc}
Z_1 & -Z_0 & Z_3 & Z_2\\
\\
Z_2 & Z_3 & -Z_0 & Z_1\\
\\
Z_3 & Z_2 & Z_1 & -Z_0\\
\\
(\beta-\gamma)Z_1 & (\beta-\gamma)Z_0 & -\alpha Z_3 & \alpha Z_2\\
\\
(\gamma - \alpha) Z_2 & \beta Z_3 & (\gamma-\alpha)Z_0 & -\beta Z_1\\
\\
(\alpha-\beta) Z_3 & -\gamma Z_2 & \gamma Z_1 & (\alpha-\beta)Z_0
\end{array}
\right)
\label{eq5.9}
\end{equation}

\bigskip

\noindent One checks that it is the union of the fiber $F_\ug$ (in the generic case) 
with the above $4$ points. 
The automorphism $\sigma$ of the characteristic variety $E_\ug$ is given by definition by the equation,
\begin{equation}
N(Z)\, \sigma (Z) = 0
\label{aut}
\end{equation}
where $\sigma(Z)$ is the column vector $\sigma(Z_\mu):= M \cdot\sigma(\ug)$ (in the variables $Z_\lambda$). One checks that $\sigma(Z)$ is already determined by the equations in (\ref{aut}) corresponding to the first three lines in $N(Z)$ which are independent of $\alpha, \beta, \gamma$ (see below). Thus $\sigma$ is in fact an automorphism of $P_3(\mathbb C)$ which is the identity on the above four points and which restricts as automorphism of $F_\ug$ for each $\ug$ generic. One checks that $\sigma$ is the product of two involutions which both restrict to $E_\ug$ (for $\ug$ generic)
\begin{equation}
\sigma = I\circ I_0
\label{eq5.11}
\end{equation}
where $I$ is the involution of the end of section 3 corresponding to $u_\mu \mapsto u^{-1}_\mu$ and where $I_0$ is given by 
\begin{equation}
I_0(Z_0)=-Z_0, \, \, I_0(Z_k)=Z_k
\label{eq5.12}
\end{equation}
for $k\in \{1,2,3\}$ and which restricts obviously to $E_\ug$ in view of (\ref{eq6.9}). Both $I$ and $I_0$ are the identity on the above four points and since $I_0$ induces the 
symmetry $\varphi(z) \mapsto \varphi(-z)$ around 
$\varphi(0)=\psi(\eta/2)$ (proposition \ref{theta}) one gets the result using proposition \ref{nice}.

\noindent  The fact that $\sigma$ does not depend on the parameters $\alpha,\beta,\gamma$ plays an  important role. Explicitly we get from the first 3 equations (\ref{aut})
\begin{equation}
\sigma(Z)_\mu =\eta_{\mu\mu} (Z^3_\mu-Z_\mu\sum_{\nu\not=\mu} Z^2_\nu-2\prod_{\lambda\not=\mu} Z_\lambda)
\label{eq4.12}
\end{equation}
for $\mu\in \{0,1,2,3\}$, where $\eta_{00}=1$ and $\eta_{nn}=-1$ for $n\in\{1,2,3\}$. $\square$

\section{Central Quadratic Forms and Generalised Cross-Products}

\noindent Let $\cala=A(V,R)=T(V)/(R)$  be a  quadratic algebra.
Its geometric data $\{E\,,\,  \sigma\,,\,\call\}$  
 is defined in
such a way that $\cala$ maps homomorphically to a cross-product algebra
obtained from  
sections of  powers of the line bundle $\call$ on powers of  the correspondence $\sigma$
(\cite{art-tat-vdb:1990}).

\noindent We shall begin by a purely algebraic result
which considerably refines the above homomorphism and lands
in a richer cross-product. We use the notations of section 4
for general quadratic algebras.

\begin{definition} \label{cent}
Let $Q \in S^2(V)$ be a symmetric bilinear form on $V^\ast$
and $C$ a component of $E \times E$. We shall 
say that $Q$ is \underline{central} on $C$ iff for all ($Z,\,Z'$) in $C$ 
and  $\omega\in R$ one has,
$$
\omega(Z,Z')\, Q(\sigma(Z'),\sigma^{-1}(Z))+Q(Z,Z')\, \omega(\sigma(Z'),\sigma^{-1}(Z))
=0
$$
\end{definition}

\noindent Let $C$ be a component of $E \times E$
globally invariant under the map
\begin{equation}
\tilde{\sigma}(Z,Z'):=\,(\sigma(Z),\sigma^{-1}(Z'))
\label{siginv}
\end{equation}
 Given a quadratic form  $Q$ central and not identically zero
on the component $C$, we define as follows
an algebra $C_Q$ as a generalised cross-product of the ring of meromorphic
functions on $C$
by the transformation $\tilde{\sigma}$.
Let $L$, $L'\in V$ be such that $L(Z)\,L'(Z')$
does not vanish identically on $C$.
We adjoin two generators $W_L$ and $W'_{L'}$ which besides the usual cross-product rules,
\begin{equation}
W_L \,f = (f\circ \tilde{\sigma}) \;W_L \,, \quad W'_{L'} \,f = (f\circ \tilde{\sigma}^{-1}) \;W'_{L'}\,, \quad \forall f\in C
\label{cropro}
\end{equation}
fulfill the following relations,
\begin{equation}
W_L \, W'_{L'}:=\pi(Z,Z')\, , \qquad W'_{L'}\, W_L := \pi(\sigma^{-1}(Z),\sigma(Z'))
\label{crossed}
\end{equation}
where the function $\pi(Z,Z')$ is given by the ratio, 
\begin{equation}
\pi(Z,Z'):=\frac{L(Z)\,L'(Z')}{Q(Z,Z')}
\label{ratio}
\end{equation}
The a priori dependence on $L$, $L'$ is eliminated by the rules,
\begin{equation}
W_{L_2}:=\frac{L_2(Z)}{L_1(Z)}\,W_{L_1}  \qquad 
W'_{L'_2}:=\,W'_{L'_1} \, \frac{L'_2(Z')}{L'_1(Z')}
\label{compare}
\end{equation}
which allow to adjoin all $W_L$ and $W'_{L'}$ for 
$L$ and $L'$ not identically zero
on the projections of $C$, without changing the algebra
and provides an intrinsic definition of $C_Q$.

\noindent Our first result is

\begin{lemma}\label{alg0}
Let $Q$ be  central and not identically zero
on the component $C$.

\noindent (i) The following equality defines a homomorphism $\rho$:
$\cala \mapsto C_Q$
$$
\sqrt{2}\;\rho(Y) :=  \frac{Y(Z)}{L(Z)}\,W_L+ W'_{L'}\,\frac{Y(Z')}{L'(Z')}\, , \qquad \forall Y \in V
$$
\noindent (ii) If $\sigma^4 \neq \bbbone$, then $\rho(Q)=1$ where 
$Q$ is viewed as an element of $T(V)/(R)$.
\end{lemma}

\noindent Formula (i) is independent of $L$, $L'$ using 
 (\ref{compare}) and reduces to $W_Y+ W'_{Y}$
when $Y$ is non-trivial on the two projections of $C$.
 It is enough to check that the  
$\rho(Y)\in C_{Q}$ fulfill the quadratic relations 
 $\omega\in R$. We view $\omega\in R$ as 
a bilinear form on  $V^\ast$. The vanishing
of the terms in $W^2$ and in $W^{'2}$ is automatic
by construction of the characteristic variety. The vanishing
of the sum of terms in $W \,W', W'\,W$ 
follows from definition \ref{cent}.

\noindent Let $\cala_\ug= C_{\mathrm{alg}}(\mathbb C^4(\ug))$ 
 at generic $\ug$, then the center of $\cala_\ug$ is generated
by the three linearly dependent quadratic elements
\begin{equation} \label{qk}
Q_m\,:=\,J_{k \ell}\,(Y_0^2 \, +\, Y_m^2)\,+ \, Y_k^2 \, - \,Y_\ell^2
\end{equation} 
with the notations of theorem \ref{iden}.

\begin{proposition} \label{central}
Let $\cala_\ug= C_{\mathrm{alg}}(\mathbb C^4(\ug))$  at generic $\ug$, then 
each $Q_m$ is central on $F_\ug \times F_\ug$ ($\subset E_\ug \times E_\ug$).
\end{proposition}

\noindent One uses (\ref{eq5.11}) to check the algebraic identity. Together
with lemma \ref{alg0}
this yields non trivial homomorphisms of $\cala_\ug$ whose
unitarity will be analysed in the next section. 
Note that for a general quadratic algebra 
$\cala=A(V,R)=T(V)/(R)$ and a quadratic form 
$Q \in S^2(V)$, such that $Q \in $ Center($\cala$),
it does not automatically follow that $Q$ is 
central on $E \times E$. For instance Proposition \ref{central}
no longer holds on $F_\ug \times \{e_\nu\}$ where $e_\nu$ is any of the 
four points of $ E_\ug$ not in $F_\ug$.

\section{Positive Central Quadratic Forms on Quad\-ratic $\ast$-Algebras}

The algebra $\cala_\ug$, $\ug \in T$ is
by construction a {\sl quadratic $\ast$-algebra} i.e. a quadratic 
complex algebra $\cala=A(V,R)$ which is 
also a $\ast$-algebra with involution $x\mapsto x^\ast$ preserving the subspace $V$ of generators. Equivalently one can take the generators of $\cala$ (spanning $V$) to be hermitian elements of $\cala$. In such a case the complex finite-dimensional vector space $V$ has a real structure given by the antilinear involution $v\mapsto j( v)$ obtained by restriction of $x\mapsto x^\ast$. Since one has $(xy)^\ast=y^\ast x^\ast$ for $x,y\in \cala$, it follows that the set $R$ of relations satisfies
\begin{equation}
(j \otimes j)( R)=t(R)
\label{eq5.3}
\end{equation}
 in $V\otimes V$ where $t:V\otimes V\rightarrow V\otimes V$ is the transposition $v\otimes w \mapsto t(v\otimes w)=w\otimes v$. This implies

\begin{lemma}\label{conj}
 The characteristic variety is stable under the 
involution $Z\mapsto j( Z)$  and one has
$$
\sigma( j(Z)) =\, j( \sigma^{-1}(Z))
$$
\end{lemma}

\noindent We let $C$ be as above an invariant component of $E \times E$
we say that $C$ is $j$-real when it is globally invariant
under the involution 
\begin{equation}
\tilde j(Z,\,Z'):=( j( Z'),\, j( Z))
\end{equation}
By lemma \ref{conj} this involution commutes with the 
automorphism $\tilde{\sigma}$ (\ref{siginv}) 
and the following turns the 
 cross-product $C_Q$ into a $\ast$-algebra,
\begin{equation}
 f^\ast(Z, \,Z'):=\overline{ f(\tilde j(Z,\,Z'))
}\,,\qquad (W_L)^\ast = W'_{j( L)}\,,\qquad (W'_{L'})^{\ast} =W_{j( L')}
\label{invol}
\end{equation}
provided that 
 $Q \in S^2(V)$ fulfills $Q=Q^\ast$.
We use the transpose of $j$, so that
\begin{equation}
 j(L)(Z)=\overline{L(j( Z))}\,, \qquad \forall Z \in V^\ast.
\label{inv2}
\end{equation}
The homomorphism
$\rho$ of lemma \ref{alg0} is a $\ast$-homomorphism. 
Composing $\rho$ with the restriction to the subset
$K= \{(Z,Z')\in C \,\vert\, Z'= j( Z)\}$ one obtains a $\ast$-homomorphism
$\theta$ of $\cala=A(V,R)$ to a twisted cross-product $C^\ast$-algebra,
$C(K) \times_{\sigma,\,\call} \mathbb {Z} $.
 Given a compact space $K$, an homeomorphism $\sigma$ of $K$ and a hermitian line bundle $\call$ on $K$ we define the  $C^\ast$-algebra $C(K) \times_{\sigma,\,\call} \mathbb {Z} $ as the twisted cross-product of $C(K)$ by the Hilbert $C^*$-bimodule associated to $\call$ and $\sigma$ (\cite{aba-eil-exe:1998}, \cite{pim:1997}).
We let for each $n \geq 0$, $\call^{\sigma^n}$ be the
hermitian line bundle pullback of $\call$
by $\sigma^n$ and (cf. \cite{art-tat-vdb:1990}, \cite{smi-sta:1992})
\begin{equation}
\call_n := \call \otimes \call^{\sigma} \otimes \cdots \otimes \call^{\sigma^{n-1}}
\label{gene2}
\end{equation}
We first define a
$\ast$-algebra as the linear span of the monomials
\begin{equation}
\xi \, W^n\, , \quad W^{\ast n} \, \eta^\ast \,,\quad \xi\,,\eta \in C(K,\call_n)
\label{gene}
\end{equation}
with product given as in (\cite{art-tat-vdb:1990}, \cite{smi-sta:1992}) for 
$(\xi_1 \, W^{n_1})\,(\xi_2 \, W^{n_2})$ so that 
\begin{equation}
(\xi_1 \, W^{n_1})\,(\xi_2 \, W^{n_2}):= (\xi_1 \otimes (\xi_2\circ{\sigma^{n_1}}) )\, W^{n_1+n_2}
\label{gene3}
\end{equation}
We use the hermitian structure of $\call_n $ to give meaning to the 
products $\eta^\ast \,\xi$ and $\xi \;\eta^\ast$ for $\xi\,,\eta \in C(K,\call_n)$.
The product then extends uniquely to an associative product of $\ast$-algebra
fulfilling the following additional rules
\begin{equation}
(W^{\ast k} \, \eta^\ast)\,( \xi \, W^k):= \, (\eta^\ast\, \xi)\circ \sigma^{-k}\,,\qquad
( \xi \, W^k)\,(W^{\ast k} \, \eta^\ast)\,:= \, \xi \;\eta^\ast
\label{gene1}
\end{equation}

\noindent The $C^\ast$-norm of $C(K) \times_{\sigma,\,\call} \mathbb {Z} $ is defined as for ordinary cross-products and due to the amenability of the group $\mathbb {Z} $
there is no distinction between the reduced and maximal norms.
The latter is obtained as the supremum of the norms in involutive
representations in Hilbert space. The natural
positive conditional expectation on the subalgebra $C(K)$ shows that the $C^\ast$-norm
restricts to the usual sup norm on $C(K)$.
 
\noindent To lighten notations
in the next statement we  abreviate $j(Z)$ as $\bar Z$,

\begin{theorem}\label{C*}
 Let $K \subset E$ be a compact
$\sigma$-invariant subset and $Q$ be central and strictly positive 
on $\{(Z,\,\bar Z);\, Z\in K\}$. Let $\call$ be the restriction to $K$
of
the dual of the tautological line bundle on $P(V^\ast)$ endowed with the unique 
hermitian metric such that $$ \displaystyle \langle L, L'\rangle =\frac{ L(Z)\,\overline{ L'(Z)}}{Q(Z,\,\bar Z)} \qquad L, L' \in V,\quad Z \in K
$$

\noindent (i) The equality $\sqrt{2}\,\theta(Y):= Y\, W + W^\ast\,\bar Y^\ast$
yields a $\ast$-homomorphism $$\theta:\cala=A(V,R)
\mapsto C(K) \times_{\sigma,\,\call} \mathbb {Z} $$ 

\noindent (ii) For any $Y \in V$
the $C^\ast$-norm of $\theta(Y)$ fulfills 
$${\rm Sup}_K \|Y\|\leq \sqrt{2}\| \,\theta(Y)\|
\leq 2\,{\rm Sup}_K \|Y\| $$

\noindent (iii) If $\sigma^4 \neq \bbbone$, then $\theta(Q)= 1$ where 
$Q$ is viewed as an element of $T(V)/(R)$.
\end{theorem}

\noindent We shall now apply this general result to the 
algebras $\cala_\ug$, 
$\ug \in T$.
We take the quadratic form
\begin{equation}
 Q(X,\,X'):=\sum X_\mu\,X'_\mu
\label{quad}
\end{equation}
in the $x$-coordinates, so that $Q$ is the canonical central element 
defining the sphere $S^3_\ug$ by the equation $Q=1$.
Proposition \ref{central} shows that $Q$ is central on $F_\ug \times F_\ug$
for generic $\ug$. The positivity of $Q$ is automatic since in the
$x$-coordinates the involution $j_\ug$
 coming from the involution of 
the quadratic $\ast$-algebra $\cala_\ug$
is simply complex conjugation
$j_\ug(Z)= \bar Z$, so that $Q(X,\,j_\ug(X))>0$
for $X \neq 0$. We thus get,

\begin{corollary} \label{II}
Let $K \subset F_\ug$ be a compact
$\sigma$-invariant subset.
The homomorphism
$\theta$ of Theorem \ref{C*}
is a unital
 $\ast$-homomorphism from $ C_{\mathrm{alg}}(S^3_\ug)$
to the cross-product $ C^{\infty}(K) \times_{\sigma,\,\call} \mathbb {Z} $. 
\end{corollary}

\noindent This applies in particular to $K=F_\ug$. It follows  that one obtains a 
non-trivial $C^\ast$-algebra $C^\ast(S^3_\ug)$
as the completion of  $ C_{\mathrm{alg}}(S^3_\ug)$
for the semi-norm,
\begin{equation}
 \| P \|:= { \rm Sup}\| \pi(P) \|
\label{norm}
\end{equation}
where $\pi$ varies through all unitary representations of 
$ C_{\mathrm{alg}}(S^3_\ug)$.

\noindent To analyse the compact
$\sigma$-invariant subsets of $F_\ug $
 for generic $\ug$, we distinguish two cases.
First note that the real curve $F_\ug \cap P_3(\mathbb {R})$
is non empty (it contains $p_0$),  
and has two connected components since it is 
invariant under the Klein group $H$ (\ref{eq6.17}).

$$  
\hbox{
\psfig{figure=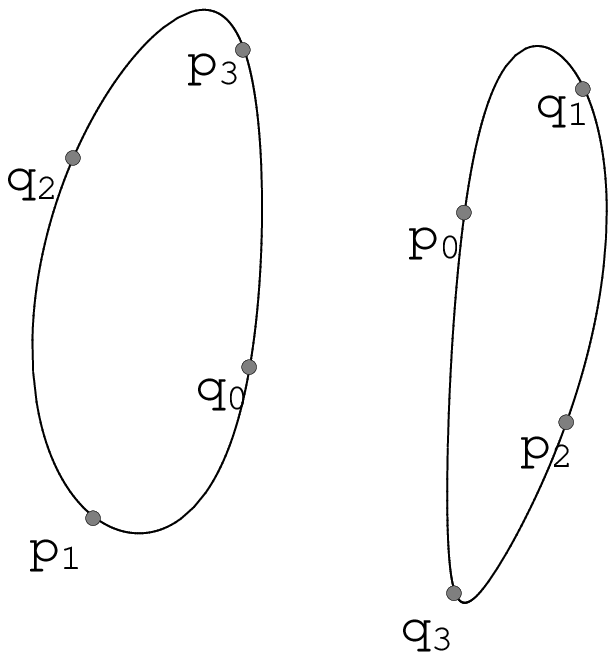}
}
$$
\centerline{Figure $2$: The Elliptic Curve $F_\ug \cap P_3(\mathbb {R})$ (odd case)}
\bigskip

\noindent We say that $\ug \in T$ is \underline{even} when $\sigma$
preserves each of the two connected components
of the real curve $F_\ug \cap P_3(\mathbb {R})$
and  \underline{odd} when it permutes them (cf. Figure $2$).
A generic $\ug \in T$ is even (cf. Figure $1$)
iff the $s_k$ of (\ref{real0})
($s_k = 1 + t_\ell\, t_m
\, , \: t_k= {\rm tan}\,\varphi_k$)
have the same sign. In that case
$\sigma$
is the square of a real translation $\kappa$
of the elliptic curve $F_\ug$
 preserving $F_\ug \cap P_3(\mathbb {R})$.

\begin{proposition} \label{siginvar} Let $\ug \in T$ be generic and even.

\noindent (i) Each connected component of $F_\ug \cap P_3(\mathbb {R})$
is a minimal compact $\sigma$-invariant subset.

\noindent (ii) Let $K \subset F_\ug$ be a compact
$\sigma$-invariant subset, then $K$ is the
sum in the elliptic curve  $F_\ug$ with origin $p_0$
of $K_{\mathbb T}=K \cap F_{\mathbb T}(\ug)^0 $ (cf. \ref{real}) with the
component $ C_\ug$ of $F_\ug \cap P_3(\mathbb {R})$
 containing $p_0$.

\noindent (iii) The cross-product $ C(F_\ug) \times_{\sigma,\,\call} \mathbb {Z} $
is isomorphic to the mapping torus of 
the automorphism $\beta
$ of the noncommutative torus ${\mathbb T}_{\eta}^2
= C_\ug \times_\sigma \mathbb {Z} $ acting on the generators by the  matrix
$\left[
\begin{array}{cc}
1& 4\\
0& 1
\end{array}
\right]
$.
\end{proposition}

\noindent More precisely with $U_j$ the generators one has
\begin{equation}
\beta(U_1):= U_1\,,\qquad \beta(U_2):= U_1^4\,U_2
\end{equation}
The mapping torus of 
the automorphism $\beta
$ is given by the algebra of continuous maps $s \in \mathbb R
\, \mapsto x(s)\in C({\mathbb T}_{\eta}^2)$
such that $ x(s+1)=\beta( x(s))\,,\; \forall s\in \mathbb R
$.

\begin{corollary} \label{H3}
Let $\ug \in T$ be generic and even, then
$F_\ug \times_{\sigma,\,\call} \mathbb {Z} $
is a noncommutative $3$-manifold with an elliptic 
 action of the three dimensional
Heisenberg Lie algebra $\frach_3$ and an invariant trace $\tau$.
\end{corollary}

\noindent We refer to \cite{rie:1989} and  \cite{aba-exe:1997} where these noncommutative manifolds were introduced and analysed in terms of crossed products by Hilbert $C^*$-bimodules. One can construct directly the action of $\frach_3$
on $C^{\infty}(F_\ug) \times_{\sigma,\,\call} \mathbb {Z}$
by choosing a constant (translation invariant) curvature connection $\nabla$, 
compatible with the metric, on the hermitian line bundle
$\call$ on $F_\ug $ (viewed in the $C^{\infty}$-category not in the
holomorphic one). The two covariant differentials $\nabla_j$
corresponding to the two vector fields $X_j$ on $F_\ug $
generating the translations of the elliptic curve, give
a natural extension of $X_j$ as the unique derivations $\delta_j$ of 
$C^{\infty}(F_\ug) \times_{\sigma,\,\call} \mathbb {Z}$
fulfilling the rules,
\begin{eqnarray}
\delta_j(f)& = & X_j(f)\,,\quad \forall f \in C^{\infty}(F_\ug)\nonumber\\
\delta_j(\xi \, W)& = & \nabla_j(\xi) \, W\,,\quad \forall \xi \in C^{\infty}(F_\ug\,,\call)
\end{eqnarray}
We let  $\delta$ be the unique derivation of 
$C^{\infty}(F_\ug) \times_{\sigma,\,\call} \mathbb {Z}$
corresponding to the grading by powers of $W$.
 It vanishes on $C^{\infty}(F_\ug)$ and 
fulfills  
\begin{equation} \label{delta}
\delta (\xi \, W^k)= i\,k\,\xi \, W^k \qquad 
\delta(W^{\ast k} \, \eta^\ast)= -i\,k\,W^{\ast k} \, \eta^\ast
 \end{equation}
Let $i\,\kappa$ be the constant curvature of the connection $\nabla$, one gets
\begin{equation}
[\delta_1,\, \delta_2]=\,\kappa\, \delta \,,\quad [\delta,\, \delta_j]=\,0
\end{equation}
which provides the required action of the Lie algebra $\frach_3$ on 
$C^{\infty}(F_\ug) \times_{\sigma,\,\call} \mathbb {Z}$.

\noindent Integration on the translation invariant volume form $dv$
of $F_\ug$ gives the $\frach_3$-invariant trace $\tau$,
\begin{eqnarray}
 \label{trace}
\tau(f)& = & \int f dv\,,\quad \forall f \in C^{\infty}(F_\ug)\nonumber\\
\tau(\xi \, W^k)& = &\tau(W^{\ast k} \, \eta^\ast)\,=\,0\,,\quad \forall k\neq 0
 \end{eqnarray}
It follows in particular that the results of \cite{ac:1980} 
apply to obtain the calculus. In particular the following gives
the ``fundamental class" as a $3$-cyclic cocycle,
\begin{equation}\label{3trace}
\tau_3(a_0,\,a_1,\,a_2
,\,a_3)=\,\sum \epsilon_{ijk}\,\tau(a_0\,\delta_i(a_1)\,\delta_j(a_2)
\,\delta_k(a_3))
\end{equation}
where the $\delta_j$ are the above derivations with $\delta_3:=\delta$.

 \noindent We shall in fact describe the
same calculus in greater generality in the last section
which will be devoted to the computation of the 
Jacobian of the homomorphism $\theta$ of corollary \ref{II}.

\noindent Similar results hold in the odd case. Then $F_\ug \cap P_3(\mathbb {R})$ 
is a minimal compact $\sigma$-invariant subset, any  compact
$\sigma$-invariant subset $K \subset F_\ug$ is 
the
sum in the elliptic curve  $F_\ug$ with origin $p_0$
of $F_\ug \cap P_3(\mathbb {R})$ with 
$K_{\mathbb T}=K \cap F_{\mathbb T}(\ug)^0 $ but the latter
is automatically invariant under the subgroup $H_0 \subset H$
of order $2$ of the Klein group $H$ (\ref{eq6.17})
\begin{equation}
H_0 :=\{ h \in H \vert \, h (F_{\mathbb T}(\ug)^0 )=F_{\mathbb T}(\ug)^0\}
\end{equation}

\bigskip
\noindent  The group law in $F_\ug$ is described geometrically as 
follows. It involves the point $q_0$. The sum $z =x + y$ 
of two points $x$ and $y$ of $F_\ug$
is $z=I_0(w)$ where $w$ is the 4th point of intersection of 
 $F_\ug$ with the plane determined by the three points $\{q_0, x, y\}$.
It commutes by construction with complex conjugation so that
$\overline{x+y}=\, \overline{ x}\,+\, \bar y \,,\quad \forall x,\,y \in F_\ug$.

\noindent  By lemma \ref{conj} the translation $\sigma$ is
imaginary for the canonical  involution $j_\ug$.
In terms of the coordinates $Z_\mu$ this involution
is  described as follows, using (\ref{change0}) (multiplied by 
$e^{i (\pi/4 -\varphi_1-\varphi_2-\varphi_3)}2^{-3/2}$)
to change variables.
Among the $3$ real numbers
\begin{eqnarray}
\lambda_k & = &\, \cos\varphi_\ell \, \cos\varphi_m\,\sin(\varphi_\ell-\varphi_m)\,,\qquad k\in \{1,2,3\}\nonumber
\end{eqnarray}
 two  have the same sign $\epsilon$ and one, $\lambda_k$,
$k\in \{1,2,3\}$, the opposite sign. Then  
\begin{equation}
 j_\ug=\,\epsilon\,I_k \,\circ \, c
\label{involE}
\end{equation}
where $c$ is  complex conjugation on the real elliptic curve $F_\ug$
(section 3)
and  $I_\mu$  the involution 
\begin{equation}
I_\mu(Z_\mu)=-Z_\mu, \, \, I_\mu(Z_\nu)=Z_\nu\,,\quad \nu \neq \mu
\label{imu}
\end{equation}
The index $k$ and the sign $\epsilon$ remain
constant when $\ug$ varies in each of the four components
of the complement of the four points $q_\mu$ in $F_{\mathbb T}(\ug)$.
The sign $\epsilon$  matters for the action of $j_\ug$ on 
linear forms as in (\ref{inv2}), but is irrelevant for the action
 on $F_\ug$. Each involution $I_\mu$ is a symmetry
$z \mapsto p-z$ in the elliptic curve $F_\ug$ and
the products $I_\mu \circ I_\nu$ form the Klein subgroup $H$ (\ref{eq6.17})
acting by translations of order two on  $F_\ug$.

 \noindent The quadratic form $Q$ of (\ref{quad})
 is given in the new coordinates by,
\begin{equation} \label{q13}
Q=\;(\prod \cos^2\varphi_{\ell}) \, \sum \, t_k\, s_k \,\,Q_k \, 
\end{equation}
 with $s_k := 1 + t_\ell\, t_m
\, , \: t_k:= {\rm tan}\,\varphi_k$
 and $Q_k$ defined by (\ref{qk}).

\noindent  Let $\ug \in T$ be generic and even
and  $v \in F_{\mathbb T}(\ug)^0$. Let $K(v)= \,v +  C_\ug $ be
 the minimal compact $\sigma$-invariant subset containing $v$
(Proposition \ref{siginvar} (ii)).
By Corollary \ref{II} we get a homomorphism,
\begin{equation} \label{ncu}
\theta_v \;:\; C_{\mathrm{alg}}(S^3_\ug) \mapsto
C^\infty({\mathbb T}_{\eta}^2)
\end{equation}
whose non-triviality will be proved below in corollary \ref{nct}.
We shall first show (Theorem \ref{gal}) that it transits through
the cross-product of the field $K_q$ of meromorphic functions on the 
elliptic curve by the subgroup of its Galois group Aut$_{\mathbb C}(K_q)$
generated by the translation $\sigma$.
 
\noindent  
For $Z = \,v + z\,$, $ z \in C_\ug$, one has using (\ref{involE}) and
(\ref{real}),
 \begin{equation}
 j_\ug(Z)
= \,I_\mu(Z-v)-\,I(v)
\label{invol2}
\end{equation}
which is a rational function $r(v,\,Z)$.
Fixing $\ug,\,v$ and substituting $Z$
and $Z'=r(v,\,Z)$ in the  formulas (\ref{crossed})
and (\ref{ratio})
of lemma \ref{alg0} with $L$ real such that  $0 \notin L(K(v))$, $L'=\,\epsilon\,L \circ I_\mu$
and $Q$ given by (\ref{q13}) 
we obtain rational formulas for a homomorphism $\tilde{\theta_v}$ of 
$ C_{\mathrm{alg}}(S^3_\ug)$ to the generalised cross-product of the field $K_q$
of meromorphic functions $f(Z)$ on the elliptic curve $F_\ug$
by $\sigma$. 
The generalised  cross-product rule
(\ref{crossed}) is given by
$W \, W':= \gamma(Z)$ where $\gamma$ is 
a rational function. Similarly
$W' \, W:=\gamma(\sigma^{-1}(Z))$.
Using integration on the cycle $K(v)$ to 
obtain a trace, together with corollary \ref{II}, we  get,

\begin{theorem} \label{gal}
 The homomorphism $\theta_v : C_{\mathrm{alg}}(S^3_\ug) \mapsto
C^\infty({\mathbb T}_{\eta}^2)$ factorises with
a homomorphism $\tilde{\theta_v}: C_{\mathrm{alg}}(S^3_\ug)
\mapsto K_q \times_\sigma \mathbb {Z}$
 to  the  generalised cross-product of the field $K_q$
of meromorphic functions on the elliptic curve $F_\ug$
 by the subgroup of the Galois group Aut$_{\mathbb C}(K_q)$
generated by $\sigma$. Its image generates the hyperfinite 
factor of type II$_1$ after weak closure relative to the trace
given by integration on the cycle $K(v)$.
\end{theorem}

\noindent Elements of $K_q$
with poles on $K(v)$ 
are unbounded and give
elements of the regular ring of affiliated operators, but all elements
of $\theta_v ( C_{\mathrm{alg}}(S^3_\ug))$ are regular on $K(v)$.
The above generalisation of the cross-product rules
(\ref{crossed}) with the rational formula for 
$W \, W':= \gamma(Z)$  is similar to the introduction 
of $2$-cocycles in the standard Brauer theory
of central simple algebras.

\section{The Jacobian of the Covering of $ S^3_\ug$}

\noindent In this section we shall analyse the morphism of $\ast$-algebras
\begin{equation}
\theta \,:\,
C_{\mathrm{alg}}(S^3_\ug)\mapsto C^\infty(F_\ug \times_{\sigma,\,\call} \mathbb {Z})
\end{equation}
of Corollary \ref{II}, by computing its Jacobian in the sense of
noncommutative differential geometry (\cite{ac:1982}).

\noindent The usual Jacobian  of a smooth map
$\varphi :\, M \mapsto N$ of manifolds is obtained as the
 ratio $\displaystyle \varphi^\ast(\,\omega_N)/\omega_M$ of the 
pullback of the volume form $\omega_N$ of the target manifold
$N$ with the volume form $\omega_M$ of the source manifold
$M$.
In noncommutative  geometry, 
differential forms $\omega$ of degree $d$ become Hochschild classes
$\;\tilde{\omega} \in HH_d(\cala)\,, \;\cala= C^\infty(M)$.
Moreover since one works  with the dual formulation in 
terms of algebras, the pullback $\varphi^\ast(\omega_N)$
is replaced by the pushforward $\varphi^{\,t}_\ast(\tilde{\omega}_N)$
under the corresponding transposed morphism of algebras 
$\;\varphi^{\,t}(f):= f \circ \varphi\,, \quad \forall f \in C^\infty(N)$.

\noindent  
The noncommutative sphere $S^3_\ug$ admits a canonical ``volume 
form" given by the Hochschild  $3$-cycle $\ch_{\frac{3}{2}}(U)$.
Our goal is to compute the push-forward,
\begin{equation}\label{goal}
\theta_\ast (\ch_{\frac{3}{2}}(U))\,\in HH_3(C^\infty(F_\ug \times_{\sigma,\,\call} \mathbb {Z}))
\end{equation}

\noindent  
 Let $\ug$ be generic and even. The noncommutative manifold
 $F_\ug \times_{\sigma,\,\call} \mathbb {Z}$ is, by Corollary \ref{H3},
a noncommutative
$3$-dimensional nilmanifold isomorphic to the mapping torus of an automorphism
of the noncommutative
$2$-torus $T^2_\eta$. Its Hohschild homology is easily computed using 
the corresponding result for the noncommutative torus (\cite{ac:1982}).
It admits in particular a canonical volume
form $V\in HH_3(C^\infty(F_\ug \times_{\sigma,\,\call} \mathbb {Z})) $ 
which corresponds to the natural class
in $HH_2(C^\infty(T^2_\eta) )$ (\cite{ac:1982}).
The volume form $V$ is
obtained in the cross-product 
 $F_\ug \times_{\sigma,\,\call} \mathbb {Z}$
from the
translation invariant $2$-form $dv$ on $F_\ug$.

\noindent To compare $\,\,\theta_\ast (\ch_{\frac{3}{2}}(U))$ with $V$
we shall pair it with the  $3$-dimensional Hochschild
 cocycle $\tau_h\in HH^3(C^\infty(F_\ug \times_{\sigma,\,\call} \mathbb {Z})) $
given, for any element $h$ of the center of $C^\infty(F_\ug \times_{\sigma,\,\call} \mathbb {Z})$, by
\begin{equation}\label{tauh}
\tau_h(a_0,\,a_1,\,a_2
,\,a_3)=\tau_3(h\,a_0,\,a_1,\,a_2
,\,a_3)
\end{equation}
where $\tau_3 \in HC^3(C^\infty(F_\ug \times_{\sigma,\,\call} \mathbb {Z})) $
is the fundamental class in cyclic cohomology defined by (\ref{3trace}).

\noindent By \cite{ac-mdv:2002a}, (2.13) p.549, the component $\ch_{\frac{3}{2}}(U)$ of the Chern 
character  is given by,
\begin{eqnarray}
\ch_{\frac{3}{2}}(U) = & - &\sum \epsilon_{\alpha\beta\gamma\delta}\cos(\varphi_{\alpha}-\varphi_{\beta}+\varphi_{\gamma}-\varphi_{\delta})\:\:x^\alpha\,dx^\beta\,dx^\gamma\,dx^\delta \, + \nonumber
\\&i& \, \sum \sin 2(\varphi_\mu-\varphi_\nu)\:\:x^\mu\,dx^\nu\,dx^\mu\,dx^\nu 
\end{eqnarray}
where $\varphi_0 :=0$. In terms of the $Y_\mu$ one gets, 
\begin{eqnarray} \label{beauty}
\ch_{\frac{3}{2}}(U)& = & \lambda\,\sum \delta_{\alpha\beta\gamma\delta}\,(s_{\alpha}-s_{\beta}+s_{\gamma}-s_{\delta})\:\:Y_\alpha\,dY_\beta\,dY_\gamma\,dY_\delta \, + \nonumber
\\&& \,\lambda\, \sum \epsilon_{\alpha\beta\gamma\delta}\;(s_\alpha\,-s_\beta)\:\:\,Y_\gamma\,dY_\delta\,dY_\gamma\,dY_\delta 
\end{eqnarray}
where $s_0:=0$, $s_k := 1 + t_\ell\, t_m
\, , \: t_k:= {\rm tan}\,\varphi_k$ and 
\begin{equation}
\delta_{\alpha\beta\gamma\delta}
=\epsilon_{\alpha\beta\gamma\delta}\,(n_{\alpha}-n_{\beta}+n_{\gamma}-n_{\delta})
\end{equation}
with $n_0=0$ and $n_k=1$. The normalization factor is 
\begin{equation}
\lambda =\,-i \;\prod \cos^2(\varphi_{k}) \, \sin(\varphi_{\ell}-\varphi_{m})
\end{equation}
Formula (\ref{beauty}) shows that, up to normalization, $\ch_{\frac{3}{2}}(U)
$ only depends 
on the fiber $F_\ug$ of $\ug$. 

\noindent Let $\ug \in T$ be generic, we assume for simplicity that $\ug$
is even, there is a similar formula
in the odd case. In our case the involutions $I$
and $I_0$ are conjugate by a real translation $\kappa$
of the elliptic curve $F_\ug$
and we let $F_\ug(0)$ be one of the two connected components of,
\begin{equation}\label{hh}
\{ Z \in F_\ug \,\vert \,I_0(Z)= \bar Z \}
\end{equation}
By Proposition \ref{siginvar} we can identify the center of 
$C^\infty(F_\ug \times_{\sigma,\,\call} \mathbb {Z}) $
with $C^\infty(F_\ug(0))$.
We assume for simplicity that $\varphi_j \in [0,\frac{\pi}{2}]$ are in
cyclic order $\varphi_k < \varphi_l <\varphi_m$ for some
 $k \in\{1,2,3\}$.

\begin{theorem} \label{vol}  
Let $h \in$ Center $(C^\infty(F_\ug \times_{\sigma,\,\call} \mathbb {Z}) )\sim
C^\infty(F_\ug(0))$. Then
$$
\langle \ch_{\frac{3}{2}}(U), \tau_h \rangle=\,6 \,\pi\, \Omega\,\int_{F_\ug(0)}h(Z)\, dR(Z)
$$
where $\Omega$ is the period given by the  elliptic integral of the first kind,
$$
\Omega=\,\int_{C_\ug} \, \frac{Z_0 dZ_k- Z_k dZ_0}{Z_\ell Z_m}
$$
and  $R$ the rational fraction,
$$
R(Z)=\, t_k\;\frac{Z_m^2}{Z_m^2+  \,c_k\,Z_l^2}
$$
 with $c_k={\rm tg}(\varphi_l)\;{\rm cot}(\varphi_k-\varphi_\ell)$.
\end{theorem}

\noindent We first obtained this result by  direct computation
using the explicit formula (\ref{3trace}) and the natural 
constant curvature connection on $\call$ given by the parametrisation
of Proposition \ref{theta} in terms of $\theta$-functions.
The pairing was first expressed 
in terms of elliptic functions and modular forms, and
the conceptual  understanding of its simplicity is at the origin of many
of the notions developped in the present paper and in particular
of the ``rational" formulation of the calculus which will be obtained
in the last section.
The geometric meaning of Theorem \ref{vol} is the computation
of the Jacobian in the sense of noncommutative geometry
of the  morphism $\theta$ as  explained above.
 The differential form 
\begin{equation}
 \omega:=\frac{Z_0 dZ_k- Z_k dZ_0}{s_k\,Z_\ell Z_m}
\end{equation}
is independent of $k$ and is, up to scale, the only holomorphic
form of type $(1,0)$ on $F_\ug$, it is invariant under the 
translations of the elliptic curve.
The integral $\Omega$ is (up to a trivial normalization factor) 
a standard elliptic integral,
it is given by an hypergeometric function in the variable 
\begin{equation}
m:=\displaystyle\frac{s_k(s_l-s_m)}{s_l(s_k-s_m)}\, 
\end{equation}
or a modular form in terms of $q$. 

\noindent The differential of $R$ is given on $F_\ug$ by 
$dR =\, J(Z) \; \omega$ where 
\begin{equation} \label{dR}
 J(Z) =\,2\,(s_m-s_l)\,c_k\,t_k\,\frac{Z_0\,Z_1
\,Z_2
\,Z_3
}{(Z_m^2+  \,c_k\,Z_l^2)^2}
\end{equation}

\noindent The period $\Omega$ does not vanish and $J(Z),\,Z \in F_\ug(0)$,
only vanishes on the $4$ ``ramification points"
necessarily present due to the symmetries.

\begin{corollary} \label{nct0}
The Jacobian of the map $\theta^{\,t}$ is given by the equality
$$
 \theta_\ast(\ch_{\frac{3}{2}}(U))= 3 \,\Omega\,J\,V 
$$
where $J$ is the element of the center $C^\infty(F_\ug(0))$ of 
$C^\infty(F_\ug \times_{\sigma,\,\call} \mathbb {Z})$
given by formula (\ref{dR}). 
\end{corollary}

\noindent This statement assumes that $\ug$
is generic in the measure theoretic sense so that $\eta$
admits good diophantine approximation (\cite{ac:1982}).
It justifies in particular the 
terminology of ``ramified covering"
applied to $\theta^{\,t}$. The function $J$ has only $4$ zeros 
on $F_\ug(0)$ which correspond to the ramification.

\noindent As shown by Theorem \ref{iden}
 the algebra $\cala_\ug$ is defined over $\mathbb R$,
i.e. admits a  natural antilinear automorphism of period two, $\gamma$
uniquely defined by
\begin{equation}
\gamma(Y_\mu):= Y_\mu\,,\quad \forall \mu
 \end{equation}
Theorem \ref{iden} also shows that $\sigma$ is defined over 
$\mathbb R$ and hence commutes with complex conjugation
$c(Z)= \bar Z$. This gives a natural real structure $\gamma$
on the algebra $C_Q$ with $C=F_\ug \times F_\ug$ and $Q$ as above,
$$
\gamma(f(Z, Z')):= \overline{f(c( Z),c( Z'))}\,,\quad 
  \gamma(W_L):=W_{c(L)} \,,\quad \gamma(W'_{L'}):=W'_{c(L')}
$$
One checks that the morphism $\rho$ of lemma \ref{alg0}
is ``real" i.e. that,
\begin{equation}
\gamma \circ \rho =\, \rho \circ \gamma
\end{equation}
Since $c(Z)= \bar Z$ reverses the orientation of $F_\ug$, 
while $\gamma$ preserves the orientation of $S^3_\ug$
it follows that
$J(\bar Z)=-J(Z)$ and $J$ necessarily vanishes on $F_\ug(0) \cap P_3(\mathbb {R})$.

\noindent Note also that for general $h$ one has
$\langle\ch_{\frac{3}{2}}(U), \tau_h \rangle \neq 0$
which shows that both $\ch_{\frac{3}{2}}(U) \in HH_3$ and
$ \tau_h \in HH^3$ are non trivial Hochschild classes. These results
 hold in the smooth algebra 
$C^\infty(S^3_\ug)$ containing the
closure of $C_{\mathrm{alg}}(S^3_\ug)$ under
holomorphic functional calculus in the $C^\ast$
algebra $C^\ast(S^3_\ug)$. We can also use Theorem \ref{vol} to show the 
non-triviality of the morphism 
$\theta_v : C_{\mathrm{alg}}(S^3_\ug) \mapsto
C^\infty({\mathbb T}_{\eta}^2)$ of (\ref{ncu}).

\begin{corollary} \label{nct}
The pullback of the fundamental 
class $[{\mathbb T}_{\eta}^2]$ of the noncommutative torus
by the homomorphism $\theta_v : C_{\mathrm{alg}}(S^3_\ug) \mapsto
C^\infty({\mathbb T}_{\eta}^2)$ of (\ref{ncu})
 is non zero, $\theta_v^\ast([{\mathbb T}_{\eta}^2]) \neq 0 \in HH^2$
provided $v$ is not a ramification point.
\end{corollary}

\noindent We have shown above the non-triviality of the Hochschild homology
and cohomology groups $HH_3(C^\infty(S^3_\ug))$ and $HH^3(C^\infty(S^3_\ug))$
by exhibiting specific elements with non-zero pairing. 
Combining the ramified cover $\pi\,=\,\theta^{\,t}$  with the natural
spectral geometry (spectral triple) on the noncommutative
$3$-dimensional nilmanifold $F_\ug \times_{\sigma,\,\call} \mathbb {Z}$ yields
a natural spectral triple on $S^3_\ug$ in the generic case.
It will be analysed in Part II, together with 
the $C^\ast$-algebra $C^\ast(S^3_\ug)$, the vanishing
of the primary class of $U$ in $K_1$, and 
the  cyclic
cohomology of $C^\infty(S^3_\ug)$.

\section{Calculus and Cyclic Cohomology}

\noindent Theorem \ref{vol} suggests the existence of a ``rational" form
of the calculus explaining the appearance of the elliptic period
$\Omega$ and the rationality of $R$. We shall show in this last section that this 
indeed the case.

\noindent Let us first go back to the general framework of twisted cross
products of the form 
\begin{equation}
\cala=C^\infty(M)\times_{\sigma,\,\call} \mathbb {Z}
\end{equation}
where $\sigma$ is a diffeomorphism of the manifold $M$.
We shall follow \cite{ac:1986b} to construct
cyclic cohomology classes from cocycles in the bicomplex
of group cohomology (with group $\Gamma=\mathbb Z$)
with coefficients in de Rham currents on $M$. The 
twist by the line bundle $\call$ introduces a non-trivial
interesting nuance.

\noindent Let $\Omega(M)$ be the
algebra  of smooth differential forms  on $M$,
endowed with the action of $\mathbb Z$
\begin{equation}
\alpha_{1,k}(\omega):=\sigma^{\ast k}\omega \,,\quad \, k\in \mathbb Z
\end{equation}
As in \cite{ac:1994} p. 219 we let $\tilde{\Omega}(M)$  
be the graded algebra obtained as the (graded) 
tensor product of $\Omega(M)$ by the exterior 
algebra $\wedge(\mathbb C [\mathbb Z]')$ on the 
  augmentation ideal $\mathbb C [\mathbb Z]'$
in the group ring $\mathbb C [\mathbb Z]$.
With $[n], n\in \mathbb Z$ the canonical basis of $\mathbb C [\mathbb Z]$, the 
augmentation $\epsilon :\mathbb C [\mathbb Z]\mapsto \mathbb C$ fulfills
$\epsilon([n])=1, \forall n\,$,  and 
\begin{equation}
\delta_n := [n]- [0]\,,\quad \, n\in \mathbb Z\,,\quad n\neq 0
\end{equation}
is a linear basis of $\mathbb C [\mathbb Z]'$. The left regular
representation of $\mathbb Z$ on $\mathbb C [\mathbb Z]$
restricts to $\mathbb C [\mathbb Z]'$
and is given on the above basis by
\begin{equation}
\alpha_{2,k}(\delta_n):= \delta_{n+k}\,-\delta_k  \,,\quad \, k\in \mathbb Z
\end{equation}
It extends to an action $\alpha_2$ of $\mathbb Z$
by automorphisms of $\wedge\mathbb C [\mathbb Z]' $.
We let $\alpha= \alpha_1 \otimes \alpha_2$ be the 
tensor product action of $\mathbb Z$
 on $\tilde{\Omega}(M)=\Omega(M)\otimes \wedge\mathbb C [\mathbb Z]' $.

\noindent We now use the hermitian line bundle $\call$ to form the 
twisted cross-product
\begin{equation}
\calc := \tilde{\Omega}(M)\times_{\alpha\,,\,\call} \mathbb {Z}
\end{equation}
We let $\call_n$ be as in (\ref{gene2}) for $n>0$ and extend its definition
for $n<0$ so that $\call_{-n}$ is the
pullback by $\sigma^n$ of the dual $\hat{\call}_n$ of $\call_n$ for all $n$.
The hermitian structure gives an antilinear isomorphism $\ast :\call_n \mapsto
\hat{\call}_n$.
The algebra $\calc$ is the 
linear span of monomials $\xi \, W^n$ where 
\begin{equation}
\xi \in C^\infty(M,\call_n) \otimes_{C^\infty(M)} \tilde{\Omega}(M)
\end{equation}
with  the product rules (\ref{gene3}), (\ref{gene1}).

\noindent Let $\nabla$ be a hermitian connection  on $\call$. We shall
turn $\calc$ into a differential graded algebra. By functoriality
$\nabla$ gives a hermitian connection on the $\call_k$ and hence
a graded derivation 
\begin{equation}\label{nabla1}
\nabla_n  : C^\infty(M,\call_n) \otimes_{C^\infty(M)} \Omega(M)\mapsto
 C^\infty(M,\call_n) \otimes_{C^\infty(M)} \Omega(M)
\end{equation}
whose square $\nabla_n^2$
is multiplication by the curvature $\kappa_n\in \Omega^2(M)$ of $\call_n$,
\begin{equation}\label{curvcoc}
\kappa_{n+m}  =  \kappa_n+\sigma^{\ast n}(\kappa_m)\,,\quad \forall n,\,m\in
\mathbb Z
\end{equation}
with $\kappa_1=\kappa\in \Omega^2(M)$ the curvature of $\call$. Ones has $d\kappa_n =0$
and we extend the differential $d$ to a graded derivation on 
$ \tilde{\Omega}(M)$
 by 
\begin{equation}
d \delta_n  = \,\kappa_n 
\end{equation}
We can then extend $\nabla_n$ uniquely to 
\begin{equation}
\tilde{\nabla}_n  : C^\infty(M,\call_n) \otimes_{C^\infty(M)} \tilde{\Omega}(M)\mapsto
 C^\infty(M,\call_n) \otimes_{C^\infty(M)} \tilde{\Omega}(M)
\end{equation}
so that it fulfills
\begin{equation}
\tilde{\nabla}_n(\xi\, \omega)= \tilde{\nabla}_n(\xi) 
\,\omega+(-1)^{{\rm deg}(\xi)}\xi \,d\omega
\,,\quad \forall \omega \in \tilde{\Omega}(M)
\end{equation}

\begin{proposition} \label{dga}

\noindent (i) The graded derivation $d$ of $\tilde{\Omega}(M)$ extends uniquely
to a graded derivation of $\calc$ such that,
$$
d(\xi \, W^n):= (\tilde{\nabla}_n (\xi)- (-1)^{{\rm deg}(\xi)}\xi \,\delta_n)\, W^n
$$

\noindent (ii) The pair $(\calc, \,d)$ is a graded differential algebra.
\end{proposition}

\noindent To construct closed graded traces on
 this differential graded algebra we follow (\cite{ac:1986b})
and consider the double complex of group cochains
(with group $\Gamma=\mathbb Z$)
with coefficients in de Rham currents on $M$. The 
cochains $\gamma \in C^{n,m}$ are totally antisymmetric maps
from $\mathbb Z^{n+1}$ to the space $\Omega_{-m}(M)$ of de Rham 
currents of dimension $-m$, which fulfill 
$$
\gamma(k_0+k,k_1+k,k_2+k,\cdots,k_n+k)=\,\sigma_{\ast}^{- k}
\gamma(k_0,k_1,k_2,\cdots,k_n)\,,\quad \forall k,\,k_j \in
\mathbb Z
$$
Besides the coboundary $d_1$ of group cohomology, given by
$$
(d_1\gamma)(k_0,k_1,\cdots,k_{n+1})=\, \sum_0^{n+1}\,
 (-1)^{j+m}\, \gamma(k_0,k_1,\cdots,\hat{k_j},\cdots,k_{n+1})
$$
and the coboundary $d_2$ of de Rham homology,
$$
(d_2\gamma)(k_0,k_1,\cdots,k_{n})=\,b(\gamma(k_0,k_1,\cdots,k_{n}))
$$
the curvatures $\kappa_n$ generate the further coboundary $d_3$ defined
on Ker $d_1$ by,
\begin{equation}
(d_3\gamma)(k_0,\cdots,k_{n+1})=\, \sum_0^{n+1}\,
 (-1)^{j+m+1}\,\kappa_{\,k_j} \gamma(k_0,\cdots,\hat{k_j},\cdots,k_{n+1})
\end{equation}
which maps Ker $d_1 \cap C^{n,m}$ to $ C^{n+1,m+2}$. Translation invariance
follows from (\ref{curvcoc}) and 
$\varphi_{\ast}(\omega C)= \varphi^{\ast-1}(\omega )\varphi_{\ast}( C)$
for $C \in \Omega_{-m}(M)$, $\omega \in \Omega^{\ast}(M)$.

\noindent To each $\gamma \in C^{n,m}$ one associates the functional
$\tilde{\gamma}$ on $\calc$ given by,
\begin{eqnarray}
\tilde{\gamma}(\,\xi \,W^n) & = & 0\,,\quad \forall\, n\neq 0
\,, \quad\xi \in \tilde{\Omega}(M)\nonumber\\
\tilde{\gamma}(\, \omega \otimes \delta_{k_1}
\cdots \delta_{k_n})& = & \langle \omega,
\gamma(0, k_1\cdots ,k_n)\rangle \,,\quad \forall k_j \in
\mathbb Z
\end{eqnarray}
and the $(n-m+1)$ linear form on $\cala=C^\infty(M)\times_{\sigma,\,\call} \mathbb {Z}$
given by,

\begin{eqnarray}
&\,&\Phi(\gamma)(a_0,a_1,\cdots,a_{n-m}) =\,
\nonumber \\&\,& \,\lambda_{n,m}\sum_0^{n-m}(-1)^{j(n-m-j)}\,
\tilde{\gamma}(da_{j+1}\cdots da_{n-m} \, a_0 \,da_1 \cdots da_{j-1}\,da_j)\qquad
\end{eqnarray}
where $\lambda_{n,m}:=\frac{n!}{(n-m+1)!}$.

\begin{lemma}\label{cyclic}
(i) The Hochschild coboundary $b\Phi(\gamma)$ is equal to $\Phi(d_1\gamma)$.

\noindent (ii) Let $\gamma \in C^{n,m}\,\cap\,$Ker $\,d_1$. Then $\Phi(\gamma)$ is a Hochschild cocycle and 

\medskip
\noindent
$\qquad \qquad\qquad\qquad B\Phi(\gamma)= \Phi(d_2\gamma)+\frac{1}{n+1}\;\Phi(d_3\gamma)$
\end{lemma}

\noindent We shall now show how the above general
framework allows to reformulate the calculus involved in Theorem \ref{vol}
in rational terms. We let $M$ be the elliptic curve $F_\ug$
where $\ug$ is generic and even. Let then $\nabla$
be an arbitrary hermitian connection on $\call$ and
$\kappa$ its curvature. We first display a 
cocycle $\gamma =\sum \gamma_{n,m}\in \sum C^{n,m}$
which reproduces the cyclic cocycle $\tau_3$.

\begin{lemma}\label{dioph}
There exists a two form $\alpha$
on  $M=F_\ug$ 
and a multiple $\lambda \,dv$
 of the translation invariant two form $dv$
such that :
 
\medskip

 (i) 
$ \qquad \qquad
\kappa_n =\,n\,\lambda \,dv+(\sigma^{\ast n}\alpha -\alpha)
\,,\quad \forall n \in
\mathbb Z
$

\medskip

 (ii) $\quad
d_2(\gamma_{\,j})=0\,,\quad d_1(\gamma_{3})=0\,,\quad d_1(\gamma_{\,1})+\frac{1}{2}\,d_3(\gamma_{3})
=0\,,\quad B\Phi(\gamma_1)=0\,,$

\medskip

 where $\gamma_{\,1}\in  C^{1,0}$ and $\gamma_{3}\in  C^{1,-2}$
 are given by

\medskip

$
\gamma_{1}(k_0,k_1):=\frac{1}{2}\,(k_1-k_0)(\sigma^{\ast k_0}\alpha+\sigma^{\ast k_1}\alpha)
\,,\quad \gamma_{3}(k_0,k_1):= k_1-k_0
\,,\quad \forall k_j \in
\mathbb Z
$

\medskip

 (iii)
The class of the cyclic cocycle $\Phi(\gamma_1)+\Phi(\gamma_3)$ is equal to $\tau_3$.
\end{lemma}
\bigskip

\noindent We use the generic hypothesis in the measure theoretic sense
to solve the ``small denominator" problem in (i). In (ii) we identify  
differential forms $\omega \in \Omega^{\,d}$ of degree $d$ with the dual currents
of dimension $2-d$.

\noindent It is a general principle explained in \cite{ac:1982} that a cyclic cocycle $\tau$
generates a calculus whose differential graded algebra is obtained as the quotient
of the universal one by the radical of $\tau$. We shall now explicitely
describe the reduced calculus obtained from the cocycle of lemma \ref{dioph} {\it(iii)}.
 We  use as above the hermitian line bundle $\call$ to form the 
twisted cross-product
\begin{equation}
\calb := \Omega(M)\times_{\alpha\,,\,\call} \mathbb {Z}
\end{equation}
of the algebra $\Omega(M)$ of differential forms on $M$ by the diffeomorphism $\sigma$.
Instead of having to adjoin the infinite number of odd elements $\delta_n$
we just adjoin two $\chi$ and $X$ as follows. 
We let $\delta$ be the derivation of $\calb$ such that
\begin{equation}
\delta (\xi \, W^n):= i \,n\,\xi \, W^n\,,\quad \forall \xi \in
C^\infty(M,\call_n) \otimes_{C^\infty(M)} \Omega(M)
\end{equation}
We adjoin $\chi$ to $\calb$ by tensoring $\calb$ with the exterior algebra
$\wedge\{\chi\}$ generated by an element $\chi$ of degree 1, and extend
the connection $\nabla$ (\ref{nabla1})
to 
the unique graded derivation $d'$ of $\Omega'=\calb \otimes \wedge\{\chi\}$ such that,
\begin{eqnarray} \label{deltaext}
d'\,\omega & = & \nabla \omega + \chi \delta(\omega)\,,\quad \forall \omega \in \calb\nonumber\\
d'\chi & = & -\,\lambda\, dv 
\end{eqnarray}
with $\lambda\, dv $ as in lemma  \ref{dioph}.
By construction, every element of $\Omega'$ is of the form
\begin{equation}
y= b_0 + b_1 \,\chi \,,\quad b_j\in \calb
\end{equation}
 One does not yet have a graded differential
algebra since  $d^{\,'2}\neq 0$. However, with $\alpha$ as in lemma  \ref{dioph}
one has 
\begin{equation}
d^{\,'2}(x)= [\,x, \,\alpha] \,,\quad \forall x \in \Omega'=\calb \otimes \wedge\{\chi\}
\end{equation}
and one can apply lemma $9$ p.$229$ of \cite{ac:1994} to get a 
differential graded algebra by adjoining the degree $1$ element
$X:= ``d1"$ fulfilling the
rules 
\begin{eqnarray} \label{X}
X^2= -\alpha \,,\quad x \,X\, y=0 \,,\quad \forall x,y \in \Omega'
\end{eqnarray}
and defining the differential $d$ by,
\begin{eqnarray} \label{deltaext1}
d\,x & = & d{\,'}\,x + [\,X, \,x]\,,\quad \forall x \in \Omega'
\nonumber\\
dX & = & 0
\end{eqnarray}
where $[\,X, \,x]$ is the graded commutator. It follows from 
lemma $9$ p.$229$ of \cite{ac:1994} that we obtain a differential 
graded algebra $\Omega^{\ast}$, generated by $\calb$, $\xi$ and $X$.
In fact using (\ref{X}) every element of $\Omega^{\ast}$ is of the 
form 
\begin{equation}
x=x_{1,1}+x_{1,2}\,X+X\,x_{2,1}+X\,x_{2,2}X\,,\quad x_{i,j}\in \Omega'
\end{equation}
and we define the functional $\int$ on $\Omega^{\ast}$ by extending the 
ordinary integral,
\begin{equation}\label{step1}
\int  \omega := \int_M \omega  \,,\quad \forall \omega \in \Omega(M)
\end{equation}
first to $\calb := \Omega(M)\times_{\alpha\,,\,\call} \mathbb {Z}
$ by
\begin{equation}\label{step2}
\int  \,\xi \, W^n\,:=0  \,,\quad \forall n\neq 0
\end{equation}
then to $\Omega'$ by
\begin{equation}  \label{step3}
\int  (b_0 + b_1 \,\chi):=\int b_1  \,,\quad \forall b_j \in \calb
\end{equation}
and finally to $\Omega^{\ast}$ as in lemma $9$ p.$229$ of \cite{ac:1994},
\begin{equation}  \label{step4}
\int (x_{1,1}+x_{1,2}\,X+X\,x_{2,1}+X\,x_{2,2}X):=\int x_{1,1}+(-1)^{{\rm deg}(\,x_{2,2})} \int x_{2,2}\, \alpha
\end{equation}

\bigskip
\begin{theorem} \label{fund} Let $M=F_\ug$,
$\nabla$, $\alpha$ be as in lemma \ref{dioph}.

\noindent The algebra $\Omega^{\ast}$ is a differential graded algebra
containing $C^{\infty}(M)\times_{\alpha\,,\,\call} \mathbb {Z}$.

\noindent  The functional $\int$ is a closed graded trace on  $\Omega^{\ast}$.

\noindent  The character of the corresponding cycle on 
$C^{\infty}(M)\times_{\alpha\,,\,\call} \mathbb {Z}$
$$
\tau(a_0, \cdots, a_{3}):= \,\int \,a_0
 \,da_1 \cdots \, da_{3}   \,,\quad \forall a_j \in 
C^{\infty}(M)\times_{\alpha\,,\,\call} \mathbb {Z}
$$
is cohomologous to the cyclic cocycle $\tau_3$.
\end{theorem}

\noindent It is worth noticing that the above calculus fits 
with \cite{ac:1980},
\cite{mdv-mic:1994}, and \cite{mdv:2001}.

\noindent Now in our case the line bundle $\call$ is holomorphic
and we can apply Theorem \ref{fund} to its canonical  hermitian connection
$\nabla$. We take the notations of section 5, with 
 $C=F_\ug \times F_\ug$, 
and $Q$ given by (\ref{q13}). This gives a particular ``rational" form 
of the calculus which explains the rationality of the answer
in Theorem \ref{vol}.
We first extend as follows the construction of
 $C_Q$. We let $\Omega(C,Q)$ be the  generalised cross-product of the
algebra $\Omega(C)$ of meromorphic differential forms (in $dZ$ and $dZ'$) on $C$
by the transformation $\tilde{\sigma}$.
The  generators $W_L$ and $W'_{L'}$ fulfill the  cross-product rules,
\begin{equation}
W_L \,\,\omega =  \tilde{\sigma}^\ast(\omega) \;W_L \,, \qquad W'_{L'} \,\omega = ( \tilde{\sigma}^{-1})^\ast(\omega) \;W'_{L'}
\label{cropro1}
\end{equation}
while (\ref{crossed}) is unchanged.
The connection $\nabla$ is the restriction to the 
subspace $\{Z'=\bar Z\}$ of the  unique graded derivation $\nabla$  on 
$\Omega(C,Q)$ which induces the usual differential on $\Omega(C)$ and satisfies,
\begin{eqnarray} \label{rational}
\nabla W_L & = & (\,d_Z\,\log L(Z)-d_Z\,\log Q(Z,Z'))\, W_L \nonumber\\
\nabla W'_{L'} & = & W'_{L'}\,(\,d_{Z'}\,\log L'({Z'})-d_{Z'}\,\log Q(Z,Z')) 
\end{eqnarray}
where $d_Z$ and $d_{Z'}$ are the (partial) differentials 
relative to the variables $Z$ and $Z'$.
Note that one needs to check that the involved differential
forms such as $d_Z\,\log L(Z)-d_Z\,\log Q(Z,Z')$ are not only
invariant under the scaling transformations $Z \mapsto \lambda Z$
but are also \underline{basic}, i.e. have zero restriction to the
fibers of the map $\mathbb C^4 \mapsto P_3(\mathbb {C})$, in both variables
$Z$ and $Z'$.
By definition the derivation 
$\delta_\kappa\,=\,\nabla^2$ of $\Omega(C,Q)$
vanishes on $\Omega(C)$ and fulfills 
\begin{equation}
\delta_\kappa(W_L)\,=\,\kappa\,W_L \,,\qquad \delta_\kappa(W'_{L'})\,=\,-\,W'_{L'}\,\kappa
\label{kappa}
\end{equation}
where 
\begin{equation} \label{curv1}
\kappa=d_Z\,d_{Z'}\,\log Q(Z,Z')
\end{equation}
 is a basic form which when restricted to the 
subspace $\{Z'=\bar Z\}$ is  the curvature. 
  We let as above $\delta$ be the derivation of  $\Omega(C,Q)$ which vanishes on $\Omega(C)$
 and is such that $\delta W_L = i\,W_L$ and $\delta W'_{L'}=-i\,W'_{L'}$. 
We proceed exactly as above and get the graded algebras 
$\Omega'=\Omega(C,Q)\otimes \wedge\{\chi\}$ obtained
by adjoining $\chi$ and $\Omega^{\ast}$ by adjoining $X$.
We define $d'$, $d$ as in  (\ref{deltaext}) and 
(\ref{deltaext1})  and the integral $\int$ 
by integration (\ref{step1}) on the subspace $\{Z'=\bar Z\}$ followed as above by 
steps (\ref{step2}), (\ref{step3}), (\ref{step4}).\\

\begin{corollary} \label{fund1} Let $\rho$:
$C_{\mathrm{alg}}(S^3_\ug)  \mapsto C_Q$ be the morphism of lemma \ref{alg0}.
  The equality
$$
\tau_{\mathrm{alg}}(a_0, \cdots, a_{3}):= \,\int \,\rho(a_0)
 \,d^{\,'}\rho(a_1) \cdots \, d^{\,'}\rho(a_{3})
$$  
defines a $3$-dimensional Hochschild cocycle $\tau_{\mathrm{alg}}$ on $C_{\mathrm{alg}}(S^3_\ug) $.

\noindent Let $h \in$ Center $(C^\infty(F_\ug \times_{\sigma,\,\call} \mathbb {Z}) )\sim
C^\infty(F_\ug(0))$. Then
$$
\langle \ch_{\frac{3}{2}}(U), \tau_h \rangle=\langle h\,\ch_{\frac{3}{2}}(U), \tau_{\mathrm{alg}} \rangle
$$

\end{corollary}

\noindent The computation of $d'$ only involves rational fractions
in the variables $Z$, $Z'$ (\ref{rational}), and the formula
(\ref{beauty}) for  $\ch_{\frac{3}{2}}(U)$ is polynomial
in the $W_L$, $W'_{L'}$. We thus obtain an a
priori reason for the rational form of the result of Theorem \ref{vol}.
The explicit computation as well as its extension to the odd case
and the degenerate cases
will be described in part II.


\end{document}